\documentclass[12pt, reqno]{amsart}
\usepackage{ifthen,amsfonts,amsmath,amssymb,epic,eepic}
\usepackage{epsfig,color}

\makeatother

\theoremstyle{definition}

\def\fnum{equation}
\newtheorem{Thm}[\fnum]{Theorem}
\newtheorem{Cor}[\fnum]{Corollary}

\newtheorem{Lem}[\fnum]{Lemma}

\newtheorem{Rem}[\fnum]{Remark}
\newtheorem{Pro}[\fnum]{Proposition}

\newcommand{\cag}{{\nu}}

\newcommand{\cat}{{\text{Cat}}}

\newcommand{\nn}{{\bf{n}}}

\newcommand{\dist}{{\text {dist}}}
\newcommand{\tsma}{\Sigma^{\theta_1,\theta_2}_{r_3,r_4}}
\newcommand{\tsm}{\Sigma^{-\pi,3\pi}_{R/2,5R/2}}
\newcommand{\sztp}{\Sigma^{0,2\pi}}

\newcommand{\Hess}{{\text {Hess}}}
\def\ZZ{{\bold Z}}
\def\RR{{\bold R}}
\def\SS{{\bold S}}

\def\CC{{\bold C }}
\newcommand{\dv}{{\text {div}}}
\newcommand{\e}{{\text {e}}}
\newcommand{\Area}{{\text {Area}}}

\newcommand{\cP}{{\mathcal{P}}}

\newcommand{\cB}{{\mathcal{B}}}

\newcommand{\cT}{{\mathcal{T}}}
\newcommand{\ba}{\beta_3}
\newcommand{\bbeta}{\bar{\beta}}

\newcommand{\gup}{{\star}}

\newcommand{\cchi}{\tilde{C}}
\newcommand{\gampa}{\gamma^{\partial}}

\newcommand{\K}{{\text{K}}}

\newcommand{\eqr}[1]{(\ref{#1})}

\begin{document}

\title[Graphical off the axis]
{The space of embedded minimal surfaces of fixed genus in a $3$-manifold I;
Estimates off the axis for disks}

\author{Tobias H. Colding}%
\address{Courant Institute of Mathematical Sciences and MIT\\
251 Mercer Street\\ New York, NY 10012 and 77 Mass. Av, Cambridge, MA 02139}
\author{William P. Minicozzi II}%
\address{Department of Mathematics\\
Johns Hopkins University\\ 3400 N. Charles St.\\ Baltimore, MD
21218}
\thanks{The first author was partially supported by NSF Grant DMS 9803253
and an Alfred P. Sloan Research Fellowship and the second author
by NSF Grant DMS 9803144 and an Alfred P. Sloan Research
Fellowship.}


\email{colding@cims.nyu.edu and minicozz@math.jhu.edu}

\maketitle


\numberwithin{equation}{section}

\section{Introduction} \label{s:s0}

This paper is the first in a series where we attempt to give a
complete description of the space of all embedded minimal surfaces
of fixed genus in a fixed (but arbitrary) closed Riemannian
$3$-manifold.  The key for understanding such surfaces is to
understand the local structure in a ball and in particular the
structure of an embedded minimal disk in a ball in $\RR^3$
(with the flat metric). This
study is undertaken here and completed in \cite{CM6}; see also \cite{CM8},
\cite{CM9} where we have surveyed our results about embedded minimal disks.
These local
results are then applied in \cite{CM7} where we
describe the general structure of fixed genus surfaces in
$3$-manifolds.

We show here that if such
an embedded minimal disk in $\RR^3$
starts off as an almost flat multi-valued graph, then it
will remain so indefinitely.

\begin{figure}[htbp]
    \setlength{\captionindent}{20pt}
    \begin{minipage}[t]{0.5\textwidth}
    \centering\input{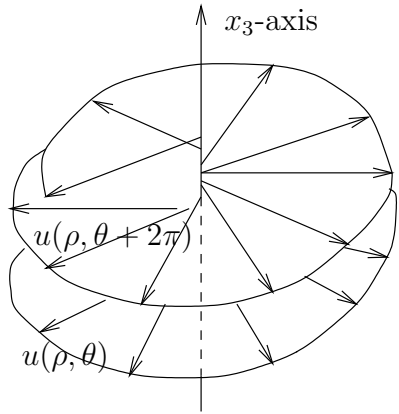}
    \caption{A multi-valued graph.}
    \label{f:1}    \end{minipage}
\end{figure}

Let $\cP$
be the universal cover of the punctured plane $\CC
\setminus \{ 0 \}$ with global (polar) coordinates $(\rho ,
\theta)$. An $N$-valued graph $\Sigma$ over the annulus $D_{r_2}
\setminus D_{r_1}$
(see fig. \ref{f:1}) is a (single-valued) graph over
\begin{equation} \label{e:defmvg}
\{ (\rho ,\theta ) \in \cP \, | \, r_1 < \rho < r_2 {\text{ and }}
|\theta| \leq \pi \, N \} \, .
\end{equation}
The middle sheet $\Sigma^M$ (an annulus with a slit as in
\cite{CM3}) is the portion over
\begin{equation}
\{ (\rho ,\theta ) \in \cP \, | \, r_1 < \rho < r_2 {\text{ and }}
0 \leq \theta \leq 2 \, \pi \} \, .
\end{equation}

\begin{Thm} \label{t:spin4ever2}
See fig. \ref{f:2}.
Given $\tau > 0$, there exist $N , \Omega, \epsilon > 0$ so:
Let $\Omega \, r_0 < 1 < R_0 / \Omega$, 
$\Sigma \subset B_{R_0}\subset \RR^3$ be
an embedded minimal disk, $\partial \Sigma \subset \partial B_{R_0}$.  If
$\Sigma$ contains an $N$-valued graph $\Sigma_g$ over
$D_1 \setminus D_{r_0}$ with gradient $\leq \epsilon$ and
$\Sigma_g \subset \{ x_3^2 \leq \epsilon^2 (x_1^2 + x_2^2) \}$, then
$\Sigma$ contains a $2$-valued graph $\Sigma_d$ over
$D_{R_0/\Omega} \setminus D_{r_0}$ with gradient $\leq \tau$ and
$(\Sigma_g)^M \subset \Sigma_d$.
\end{Thm}

\begin{figure}[htbp]
    \setlength{\captionindent}{20pt}\begin{minipage}[t]{0.5\textwidth}
\centering\input{shn2.pstex_t}
    \caption{Theorem \ref{t:spin4ever2} - extending a small multi-valued graph
    in a disk.}
\label{f:2}    \end{minipage}\begin{minipage}[t]{0.5\textwidth}
    \centering\input{shn3.pstex_t}
    \caption{Theorem \ref{t:blowupwindinga} - 
finding a small multi-valued graph in
    a disk near a point of large curvature.}
    \label{f:3}    \end{minipage}
\end{figure}

Theorem \ref{t:spin4ever2} is particularly useful when combined with a result
from \cite{CM4} asserting that an embedded
minimal disk with large curvature at a point
contains a small almost flat multi-valued graph nearby.
Namely:

\begin{Thm} \label{t:blowupwindinga}
\cite{CM4}.
See fig. \ref{f:3}.
Given $N , \omega>1$, $\epsilon > 0$, there exists
$C=C(N,\omega,\epsilon)>0$ so: Let
$0\in \Sigma^2\subset B_{R}\subset \RR^3$ be an embedded minimal
disk, $\partial \Sigma\subset \partial B_{R}$. If
$\sup_{B_{r_0} \cap \Sigma}|A|^2\leq 4\,C^2\,r_0^{-2}$
and $|A|^2(0)=C^2\,r_0^{-2}$ for some $0<r_0<R$, then there exist
$ \bar{R} < r_0 / \omega$ and (after a rotation)
an $N$-valued graph $\Sigma_g \subset \Sigma$ over $D_{\omega \bar{R} }
\setminus D_{\bar{R} }$ with gradient $\leq \epsilon$, and
$\dist_{\Sigma}(0,\Sigma_g) \leq 4 \, \bar{R}$.
\end{Thm}

Combining these two results with a standard blow up argument gives:

\begin{Thm} \label{t:blowupwinding0}
\cite{CM4}.
Given $N\in \ZZ_+$, $\epsilon > 0$, there exist
$C_1,\,C_2>0$ so: Let
$0\in \Sigma^2\subset B_{R}\subset \RR^3$ be an embedded minimal
disk, $\partial \Sigma\subset \partial B_{R}$. If
$\max_{B_{r_0} \cap \Sigma}|A|^2\geq 4\,C_1^2\,r_0^{-2}$ for some
$0<r_0<R$, then there exists
(after a rotation)
an $N$-valued graph $\Sigma_g \subset \Sigma$ over $D_{R/C_2}
\setminus D_{2r_0}$ with gradient $\leq \epsilon$
and
$\Sigma_g \subset \{ x_3^2 \leq \epsilon^2 \, (x_1^2 + x_2^2) \}$.
\end{Thm}

The multi-valued graphs given by Theorem \ref{t:blowupwinding0}
should be thought of (see \cite{CM6}) as the basic building blocks
of an embedded minimal disk. In fact, one should think of such a
disk as being built out of such graphs by stacking them on top of
each other. It will follow from Proposition \ref{l:grades1} that
the separation between the sheets in such a graph grows
sublinearly.

\begin{figure}[htbp]
    \setlength{\captionindent}{20pt}
    \begin{minipage}[t]{0.5\textwidth}
\centering\input{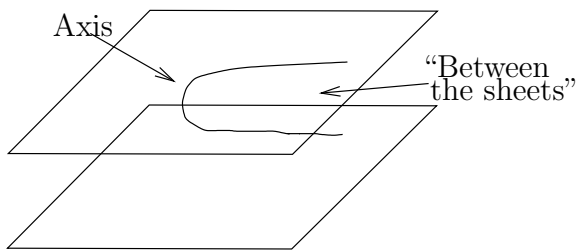}
    \caption{The estimate between the sheets: 
	Theorem \ref{t:tslab}.}
\label{f:4}    \end{minipage}
\end{figure}

An important component of the proof of
Theorem \ref{t:spin4ever2} is a version of it for stable minimal
annuli with slits that start off as multi-valued graphs.
Another component is a curvature
estimate ``between the sheets'' for embedded minimal disks in
$\RR^3$; see fig. \ref{f:4}.
We will think of an axis for such a disk $\Sigma$ as a point or
curve away from which the
surface locally (in an extrinsic ball) has more than one
component.  With this weak notion of an axis, our estimate is that if
one component of $\Sigma$ is sandwiched between two
others that connect to an axis, then the one that
is sandwiched has curvature estimates; see Theorem \ref{t:tslab}.
The example to keep in
mind is a helicoid and the components are
``consecutive sheets'' away from the axis.

Theorems \ref{t:spin4ever2}, \ref{t:blowupwindinga},
\ref{t:blowupwinding0}
are local and are for simplicity
stated and proven only in $\RR^3$ although
they can with only very minor changes easily be seen to hold for
minimal disks in a sufficiently small ball in any
given fixed Riemannian $3$-manifold.

The paper is divided into $4$ parts. In Part \ref{p:p1}, we show
the curvature estimate ``between the sheets''
when the
disk is in a thin slab. In Part \ref{p:p2}, we
will show that certain stable disks with interior
boundaries starting off as multi-valued graphs remain very flat
(cf. Theorem \ref{t:spin4ever2}). This result will be needed together
with Part \ref{p:p1} in
Part \ref{p:p3} to generalize the results of Part \ref{p:p1}
to when the disk is not anymore assumed to lie in a
slab. Part \ref{p:p2} will also be used together with Part \ref{p:p3}
in Part \ref{p:p4}
to show Theorem \ref{t:spin4ever2}.

Let $x_1 , x_2 , x_3$ be the standard coordinates on $\RR^3$ and
$\Pi : \RR^3 \to \RR^2$ orthogonal projection to $\{ x_3 = 0 \}$.
For $y \in S \subset \Sigma \subset \RR^3$ and $s > 0$, the
extrinsic and intrinsic balls and tubes are
\begin{alignat}{2}
B_s(y) &= \{ x \in \RR^3 \, | \, |x-y| < s \} \, , \, & T_s(S) &=
\{ x \in \RR^3 \, | \, \dist_{\RR^3} (x , S) < s \} \, , \\
\cB_s(y) &= \{ x \in \Sigma \, | \, \dist_{\Sigma} (x , y) < s \}
\, , \, & \cT_s (S) &= \{ x \in \Sigma \, | \, \dist_{\Sigma} (x ,
S) < s \} \, .
\end{alignat}
$D_s$ denotes the disk $B_s(0) \cap \{ x_3 = 0 \}$.
$\K_{\Sigma}$ the sectional curvature of a smooth compact surface
$\Sigma$ and when
$\Sigma$ is immersed $A_{\Sigma}$ will be its second fundamental form.
When $\Sigma$ is oriented, $\nn_{\Sigma}$ is the unit normal.
We will often consider the
intersection of curves and
surfaces with extrinsic balls.   We
assume that these intersect  transversely since this can
be achieved by an arbitrarily small perturbation of the
radius.


\setcounter{part}{0}
\numberwithin{section}{part} 
\renewcommand{\rm}{\normalshape} 
\renewcommand{\thepart}{\Roman{part}}
\setcounter{section}{1}

\part{Minimal disks in a slab} \label{p:p1}

\begin{figure}[htbp]
    \setlength{\captionindent}{20pt}
    \begin{minipage}[t]{0.5\textwidth}
\centering\input{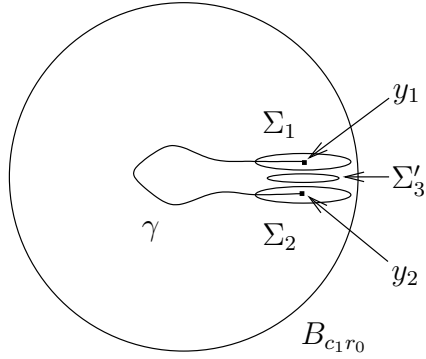}
    \caption{$y_1$, $y_2$, $\Sigma_1$, $\Sigma_2$, $\Sigma_3'$, 
and $\gamma$ in Theorem \ref{t:tslab}.}
\label{f:newfig}    \end{minipage}
\end{figure}

Let  $\gamma_{p,q}$ denote the line segment from $p$ to
$q$ and $\overline{p,q}$ the ray from $p$ through $q$. A curve
$\gamma$ is $h$-{\it{almost monotone}} if given $y \in \gamma$,
then $B_{4 \, h}(y) \cap \gamma$ has only one component which
intersects $B_{2\, h }(y)$.
Our curvature estimate ``between the sheets'' is
(see fig. \ref{f:newfig}):

\begin{Thm} \label{t:tslab}
There exist $c_1 \geq 4$, $2 c_2 < c_4 < c_3 \leq 1$ so:
Let $\Sigma^2 \subset B_{c_1 \,r_0}$ be an embedded minimal disk with
$\partial \Sigma \subset
\partial B_{c_1 \, r_0}$ and $y \in \partial B_{2\, r_0}$.
Suppose $\Sigma_1 , \Sigma_2 , \Sigma_3$ are distinct
components of $B_{ r_0}(y) \cap \Sigma$ and $\gamma \subset (
B_{r_0} \cup T_{c_2 \,r_0}(\gamma_{0,y}) ) \cap \Sigma$ is a curve
with $\partial \gamma = \{ y_1 , y_2 \}$ where $y_i \in B_{c_2 \,
r_0 }(y) \cap \Sigma_i$ and each component of $\gamma \setminus
B_{r_0}$ is $c_2 \, r_0$-almost monotone. Then any component
$\Sigma_3'$ of $B_{c_3 \, r_0 }(y) \cap \Sigma_3$ with $y_1 ,
y_2$ in distinct components of $B_{c_4 \, r_0}(y) \setminus
\Sigma'_3$ is a graph.
\end{Thm}

The idea for the proof of Theorem \ref{t:tslab}
 is to show that if this were not the case, then
we could find an embedded stable disk that would be almost
flat and lie in the complement of the original  disk. In fact, we
can choose the stable disk to be sandwiched between the two
components as well. The flatness would force the stable disk to
eventually cross the axis in the original disk, contradicting that
they were disjoint.

In this part, we prove Theorem
\ref{t:tslab} when the surface is in a slab, illustrating the key
points (the full theorem, using the results of this part, will be 
proven later).
Two simple facts about minimal surfaces in a slab will be
used: (1) Stable surfaces in a slab must be graphical away
from their boundary (see Lemma \ref{l:lone} below) and (2) the
maximum principle, and catenoid foliations in particular, force
these surfaces to intersect a narrow cylinder about every vertical
line (see the appendix).

\begin{Lem} \label{l:lone}
Let $\Gamma \subset \{ |x_3| \leq \beta \, h \}$ be a stable
embedded minimal surface. There exist $C_g , \beta_s > 0$ so
if $\beta \leq \beta_s$ and $E$ is a component of $\RR^2 \setminus
T_{h}( \Pi (\partial \Gamma ))$, then each component of
$\Pi^{-1}(E) \cap \Gamma$ is a graph over $E$ of a function $u$
with $|\nabla_{\RR^2} u| \leq C_g \, \beta$.
\end{Lem}

\begin{proof}
If $\cB_h (y) \subset \Gamma$, \cite{Sc} gives that $|A|^2 \leq
C_s \, h^{-2}$ on $\cB_{h/2}(y)$. Since $\Delta_{\Gamma} x_3 = 0$,
\cite{ChY} yields
\begin{equation} \label{e:grepi}
\sup_{\cB_{h/4}(y)} |\nabla_{\Gamma} x_3| \leq \bar{C}_g \, h^{-1}
\, \sup_{\cB_{h/2}(y)} | x_3| \leq \bar{C}_g \, \beta \, ,
\end{equation}
where $\bar{C}_g = \bar{C}_g (C_s)$. Since $|\nabla_{\RR^2} u|^2 =
|\nabla_{\Gamma} x_3|^2 \, / \, (1 - |\nabla_{\Gamma} x_3|^2)$,
this gives the lemma.
\end{proof}

The next lemma shows that if an embedded minimal disk $\Sigma$ in
the intersection of a ball with a thin slab is not graphical near
the center, then it contains a curve $\gamma$ coming close to the
center and connecting two boundary points which are close in
$\RR^3$ but not in $\Sigma$. The constant $\beta_A$ is defined in
\eqr{e:lowh}.

\begin{Lem} \label{l:la}
Let $\Sigma^2 \subset B_{60 \, h} \cap \{ |x_3| \leq \beta_A \, h
\}$ be an embedded minimal disk with $\partial \Sigma
\subset
\partial B_{60 \, h}$ and let $z_b \in \partial B_{50 \, h }$. If
a component $\Sigma'$ of $B_{5 \, h } \cap \Sigma$ is not a graph,
then there are distinct components $S_1,\, S_2$ of $B_{8 \, h
}(z_b)\cap \Sigma$, $z_i \in B_{h / 4} (z_b) \cap S_i$ and a curve
$\gamma \subset (B_{30 \, h } \cup T_{ h} (\gamma_{q,z_b}) )\cap
\Sigma$ with $\partial \gamma = \{ z_1 , z_2 \}$ and $\gamma \cap
\Sigma' \ne \emptyset$. Here $q \in B_{50\, h}(z_b)\cap \partial
B_{30\, h}$.
\end{Lem}

\begin{figure}[htbp]
    \setlength{\captionindent}{20pt}
    \begin{minipage}[t]{0.5\textwidth}
    \centering\input{shn5a.pstex_t}
    \caption{Proof of Lemma \ref{l:la}:
    Vertical plane tangent to $\Sigma$ at $z$.  Since $\Sigma$ is
    minimal, we get locally near $z$ on one side of the plane
    two different components.
    Next place a catenoid foliation centered at $y$ and tangent
    to $\Sigma$ at $z$.}
    \label{f:la1}    \end{minipage}\begin{minipage}[t]{0.5\textwidth}
\centering\input{shn5b.pstex_t}
    \caption{Proof of Lemma \ref{l:la}:  Step 1:
    Using the catenoid foliation, we build out the curve to scale $h$.}
\label{f:la2}    \end{minipage}
\end{figure}

\begin{proof}
See fig. \ref{f:la1}.
Since $\Sigma'$ is not graphical, we can find $z \in \Sigma'$ with
$\Sigma$ vertical at $z$ (i.e., $|\nabla_{\Sigma} x_3|(z) = 1$).
Fix $y \in \partial B_{4 \, h}(z)$ so $\gamma_{y,z}$ is normal to
$\Sigma$ at $z$. Then $f_y (z) = 4 \, h $ (see \eqr{e:deffy}). Let
$y'$ be given by that $y' \in \partial B_{10 \, h}(y)$ and $z \in
\gamma_{y,y'}$. The first step is to use the catenoid foliation
$f_y$ to build the desired curve on the scale of $h$; see fig. \ref{f:la2}. The second
and third steps will bring the endpoints of this curve out near
$z_b$.

Any simple closed curve $\sigma \subset \Sigma \setminus \{ f_y >
4 \, h \}$ bounds a disk $\Sigma_{\sigma} \subset \Sigma$. By
Lemma \ref{l:apb}, $f_y$ has no maxima on $\Sigma_{\sigma} \cap \{
f_y > 4 \, h \}$ so $\Sigma_{\sigma} \cap \{ f_y > 4 \, h \} =
\emptyset$. On the other hand, by Lemma \ref{l:apa}, we get a
neighborhood $U_z \subset \Sigma$ of $z$ where $U_z \cap \{ f_y =
4 \, h \} \setminus \{ z \}$ is the union of $2n \geq 4$ disjoint
embedded arcs meeting at $z$. Moreover, $U_z \setminus \{ f_y \geq
4 \, h \}$ has $n$ components $U_1 , \dots , U_n$ and
$\overline{U_i} \cap \overline{U_j} = \{ z \}$ for $i \ne j$. If a
simple curve $\tilde{\sigma}_z \subset \Sigma \setminus \{ f_y
\geq 4 \, h \}$ connects $U_1$ to $U_2$, connecting $\partial
\tilde{\sigma}_z$ by a curve in $U_z$ gives a simple closed curve
$\sigma_z \subset \Sigma \setminus \{ f_y
> 4 \, h \}$
with $\tilde{\sigma}_z \subset \sigma_z$ and $\sigma_z \cap \{ f_y
\geq 4 \,h \} = \{ z \}$. Hence, $\sigma_z$ bounds a disk
$\Sigma_{\sigma_z} \subset \Sigma \setminus \{ f_y > 4 \,h \}$. By
construction, $U_z \cap \Sigma_{\sigma_z} \setminus \cup_i
\overline{U_i} \ne \emptyset$, which is a contradiction. This
shows that $U_1 , U_2$ are contained in components $\Sigma_{4 \,
h}^1 \ne \Sigma_{4 \, h}^2$ of $\Sigma \setminus \{ f_y \geq 4 \,
h \}$ with $z \in \overline{\Sigma_{4 \, h}^1} \cap
\overline{\Sigma_{4 \, h}^2}$. For $i=1,2$, Lemma \ref{l:apb} and
\eqr{e:lowh} give $y_i^a \in B_{ h / 4 } (y) \cap \Sigma_{4 \,
h}^i$. Corollary \ref{c:c3} gives $\cag_i \subset
T_{h}(\gamma_{y,y'}) \cap \Sigma$ with $\partial \cag_i = \{ y_i^a
, y_i^b \}$ where $y_i^b \in B_{h / 4}(y')$. There are now two
cases: If $y_1^b$ and $y_2^b$ do not connect in $B_{4 \, h} (y')
\cap \Sigma$, then take $\gamma_0 \subset B_{ 5 \, h} (y) \cap
\Sigma$ from $y_1^a$ to $y_2^a$ and set $\gamma_a = \cag_1 \cup
\gamma_0 \cup \cag_2$ and $y_i = y_i^b$. Otherwise, if
$\hat{\gamma}_0 \subset B_{ 4 \, h } (y') \cap \Sigma$ connects
$y_1^b$ and $y_2^b$, set $\gamma_a = \cag_1 \cup \hat{\gamma}_0
\cup \cag_2$ and $y_i = y_i^a$. After possibly switching $y$ and
$y'$, we get a curve $\gamma_a \subset (T_{h}(\gamma_{y,y'}) \cup
B_{5 \, h}(y') ) \cap \Sigma$ with $\partial \gamma_a = \{ y_1 ,
y_2 \} \subset B_{h / 4}(y)$ and $y_i \in S_i^a$ for components
$S_1^a \ne S_2^a$ of $B_{4 \, h}(y) \cap \Sigma$. This completes
the first step.

\begin{figure}[htbp]
    \setlength{\captionindent}{20pt}
    \begin{minipage}[t]{0.5\textwidth}
\centering\input{shn5c.pstex_t}
    \caption{Proof of Lemma \ref{l:la}:  Step 2:  $y_1$ and $y_2$
cannot connect in the half-space $H$ since this would give a point in
$\Sigma_{1,2}$ far from $\partial \Sigma_{1,2}$, contradicting
Corollary \ref{c:c3}.}
    \label{f:la3}    \end{minipage}
\end{figure}

Second, we use the maximum principle to restrict the possible
curves from $y_1$ to $y_2$;  see fig. \ref{f:la3}. Set
\begin{equation}
H = \{ x \, | \, \langle y-y' , x-y \rangle > 0 \} \, .
\end{equation}
If $\eta_{1,2} \subset T_{h}(H) \cap \Sigma$ connects $y_1$ and
$y_2$, then $\eta_{1,2} \cup \gamma_a$ bounds a disk $\Sigma_{1,2}
\subset \Sigma$. Since $\eta_{1,2} \subset T_{h}(H)$, $\partial
B_{8 \, h}(y') \cap \partial \Sigma_{1,2}$ consists of an odd
number of points in each $S_i^a$ and hence $\partial B_{8 \,
h}(y') \cap \Sigma_{1,2}$ contains a curve from $S_1^a$ to
$S_2^a$. However, $S_1^a$ and $S_2^a$ are distinct components of
$B_{4 \, h}(y) \cap \Sigma$, so this curve contains
\begin{equation} \label{e:mcpo}
y_{1,2} \in \partial B_{4 \, h}(y) \cap
\partial B_{8 \, h}(y') \cap \Sigma_{1,2}
\, .
\end{equation}
By construction, $\Pi(y_{1,2})$ is in an unbounded component of
$\RR^2 \setminus T_{h/ 4} (\Pi (\partial \Sigma_{1,2}))$,
contradicting Corollary \ref{c:c2}. Hence, $y_1$ and $y_2$ cannot
be connected in $T_{h}(H) \cap \Sigma$.

Third, we extend $\gamma_a$. There are two cases: (A) If $z_b \in
H$, Corollary \ref{c:c3} gives
\begin{equation}
\tilde{\cag}_1 , \tilde{\cag}_2 \subset T_{h}(\gamma_{y,z_b}) \cap
\Sigma \subset T_{h}(H) \cap \Sigma
\end{equation}
from $y_1 , y_2$ to $ z_1 , z_2 \in B_{h / 4}(z_b),$ respectively.
(B) If $z_b \notin H$, then fix $z_c \in \partial B_{20 \,
h}(y) \cap \Pi (\partial H)$ on the same side of $\Pi
(\overline{y,y'})$ as $\Pi(z_b)$ and fix $z_d \in \partial
B_{10 \, h}(z_c) \setminus H$ with $\gamma_{z_c , z_d}$
orthogonal to $\partial H$ (so $\Pi (y'),\Pi (y) ,z_c , z_d$
form a $10\, h$ by $20 \, h$ rectangle). Corollary \ref{c:c3}
gives
\begin{equation}
\tilde{\cag}_1 , \tilde{\cag}_2 \subset T_{h}( \gamma_{y,z_c} \cup
\gamma_{z_c,z_d} \cup \gamma_{z_d,z_b}) \cap \Sigma
\end{equation}
from $y_1 , y_2$ to $ z_1 , z_2 \in B_{h / 4}(z_b)$, respectively.
In either case, set $\gamma = \tilde{\cag}_1 \cup \gamma_a \cup
\tilde{\cag}_2$. Set $q=
\partial B_{30 \, h}(y) \cap \gamma_{y,z_b}$ (in (A)) or
$q=\partial B_{30 \, h}(y) \cap \gamma_{z_c,z_b}$ (in (B)).
Applying Corollary \ref{c:c2} as above, $ z_1 , z_2$ are in
distinct components of $B_{8 \, h}(z_b) \cap \Sigma$.
\end{proof}

The next result illustrates the main ideas for Theorem
\ref{t:tslab} in the simpler case where $\Sigma$ is in a slab. Set
$\ba = \min \{ \beta_A , \beta_s , \tan \theta_0 / (2 \, C_g )
\}$; $C_g , \beta_s$ are defined in Lemma \ref{l:lone}, $\theta_0$
in \eqr{e:deft0}, and $\beta_A$ in \eqr{e:lowh}.

\begin{Pro} \label{p:pslabn1}
Let $\Sigma \subset B_{4 \, r_0} \cap \{ |x_3| \leq \ba \, h \}$
be an embedded minimal disk with $\partial \Sigma \subset
\partial B_{4\, r_0}$ and let $y \in \partial B_{2\, r_0}$.
Suppose that $\Sigma_1 , \Sigma_2 , \Sigma_3$ are distinct
components of $B_{r_0}(y) \cap \Sigma$ and $\gamma \subset (
B_{r_0} \cup T_{h}(\gamma_{0,y}) ) \cap \Sigma$ is a curve with
$\partial \gamma = \{ y_1 , y_2 \}$ where $y_i \in B_{h}(y) \cap
\Sigma_i$ and each component of $\gamma \setminus B_{r_0}$ is
$h$-almost monotone. Then any component $\Sigma_3'$ of $B_{r_0 -
80 \, h}(y) \cap \Sigma_3$ for which $y_1 , y_2$ are in distinct
components of $B_{5 \, h}(y) \setminus \Sigma'_3$ is a graph.
\end{Pro}

\begin{proof}
We will suppose that $\Sigma_3'$ is not a graph and deduce a
contradiction. Fix a vertical point $z \in \Sigma_3'$. Define $z_0
, y_0 , y_b$ on the ray $\overline{0,y}$ by $z_0 = \partial B_{3\,
r_0 - 21 \, h} \cap \overline{0,y}$, $y_0 = \partial B_{3\, r_0 -
10 \, h} \cap \overline{0,y}$, and $y_b = \partial B_{4 \, r_0}
\cap \overline{0,y}$. Set $z_b = \partial B_{50 \, h}(z) \cap
\gamma_{z, z_0}$. Define the half-space
\begin{equation}
H = \{ x \, | \, \langle x - z_0 , z_0 \rangle > 0 \} \, .
\end{equation}

The first step is to find a simple curve $\gamma_3 \subset \left(
B_{r_0 - 20 \, h}(y) \cup T_h(\gamma_{y,y_b}) \right) \cap \Sigma$
which can be connected to $\Sigma_3'$ in $B_{r_0 - 20 \, h}(y)
\cap \Sigma$, with $\partial \gamma_3 \subset \partial \Sigma$,
and so $\partial B_{ r_0 - 10 \, h}(y) \cap \gamma_3$ consists of
an odd number of points in each of two distinct components of $H
\cap \Sigma$. To do that, we begin by applying Lemma \ref{l:la} to
get $q \in B_{50 \, h}(z_b) \cap \partial B_{30 \, h }(z)$,
distinct components $S_1 ,S_2$ of $ B_{8 \, h}(z_b) \cap \Sigma$
with $z_i \in B_{h/ 4} (z_b) \cap S_i$, and a curve
\begin{equation} \label{e:onemS}
\gamma_3^{\gup} \subset (B_{30 \, h }(z) \cup T_{ h}
(\gamma_{q,z_b}) ) \cap \Sigma , \,
\partial \gamma_3^{\gup} = \{ z_1 , z_2 \} \, , \,
\gamma_3^{\gup} \cap \Sigma_3' \ne \emptyset \, .
\end{equation}
Corollary \ref{c:c3} gives $h$-almost monotone curves $\nu_{1},
\nu_{2} \subset T_h (\gamma_{z_b,z_0} \cup \gamma_{z_0, y_b} )
\cap \Sigma$ from $z_1 , z_2$, respectively, to $\partial \Sigma$.
Then $\gamma_3 = \nu_{1} \cup \gamma_3^{\gup} \cup \nu_{2}$
extends $\gamma_3^{\gup}$ to $\partial \Sigma$. Fix $z^+ \in
B_{h}(y_0) \cap \nu_{1}$ and $z^- \in B_{h}(y_0) \cap \nu_{2}$. We
will show that $z^+ , z^-$ do not connect in $H \cap \Sigma$. If
$\eta_{+}^{-} \subset H \cap \Sigma$ connects $z^+ $ and $ z^-$,
then $\eta_{+}^{-}$ together with the portion of $\gamma_3$ from
$z^+$ to $z^-$ bounds a disk $\Sigma_{+}^{-} \subset \Sigma$.
Using the almost monotonicity of each $\nu_{i}$, $\partial B_{50
\, h}(z) \cap \partial \Sigma_{+}^{-}$ consists of an odd number
of points in each $S_i$. Consequently, a curve
$\sigma_{+}^{-} \subset \partial B_{50 \, h}(z) \cap
\Sigma_{+}^{-}$ connects $S_1$ to $S_2$ and so
$\sigma_{+}^{-} \setminus B_{8 \, h}(z_b) \ne \emptyset$.  This
would contradict Corollary \ref{c:c2} and we conclude that there are
distinct components $\Sigma_{H}^+$ and $\Sigma_{H}^-$ of $H \cap
\Sigma$ with $z^{\pm} \in \Sigma_{H}^{\pm}$. Finally, removing any
loops in $\gamma_3$ (so it is simple) gives the desired
curve.

The second step is to find disjoint stable disks $\Gamma_1 ,
\Gamma_2 \subset B_{r_0 - 2 \, h}(y) \setminus \Sigma$ with
$\partial \Gamma_i \subset \partial B_{r_0 - 2 \, h}(y)$ and
graphical components $\Gamma_i'$ of $B_{r_0 - 4 \, h}(y) \cap
\Gamma_i$ so $\Sigma_3'$ is between $\Gamma_1' , \Gamma_2'$
and $y_1 , y_2 , \Sigma_3'$ are each in their own component of
$B_{r_0 - 4 \, h}(y) \setminus (\Gamma_1' \cup \Gamma_2')$. To
achieve this, we will solve two Plateau problems using $\Sigma$ as
a barrier and then use that $\Sigma_3'$ separates $y_1 , y_2$ near
$y$ to get that these are in different components. Let $\Sigma_1'
, \Sigma_2'$ be the components of $B_{r_0 - 2 \, h}(y) \cap
\Sigma$ with $y_1 \in \Sigma_1' , y_2 \in \Sigma_2'$. By the
maximum principle, each of these is a disk. Let $\Sigma_{y_2}$ be
the component of $B_{3\, h}(y_1) \cap \Sigma$ with $y_2 \in
\Sigma_{y_2}$.  Since $y_1 \notin \Sigma_{y_2}$,
Lemma \ref{l:apb} gives $y_2' \in \Sigma_{y_2} \setminus
N_{\theta_0}(y_1)$   with $\theta_0 > 0$ from \eqr{e:deft0}.
 Hence, the vector $y_1-y_2'$ is nearly
orthogonal to the slab, i.e.,
\begin{equation} \label{e:nv}
|\Pi (y_2' - y_1)| \leq |y_2' - y_1| \, \cos \theta_0 \, .
\end{equation}
Since $\Sigma_3'$ separates $y_1 , y_2$ in $B_{5\, h}(y)$, we get
$y_3 \in \gamma_{y_1,y_2'} \cap \Sigma_3'$. Fix a component
$\Omega_1$ of $B_{r_0 - 2 \, h}(y) \setminus \Sigma$ containing a
component of $\gamma_{y_1,y_3} \setminus \Sigma$ with exactly one
endpoint in $\Sigma_1'$. By \cite{MeYa}, we get a stable embedded
 disk $\Gamma_1 \subset \Omega_1$ with $\partial \Gamma_1 =
\partial \Sigma_1'$.
Similarly, let $\Omega_2$ be a component of $B_{r_0 - 2 \, h}(y)
\setminus (\Sigma \cup \Gamma_1)$ containing a component of
$\gamma_{y_3,y_2'} \setminus (\Sigma \cup \Gamma_1)$ with exactly
one endpoint in $\Sigma_2'$. Again by \cite{MeYa}, we get a stable
embedded disk $\Gamma_2 \subset \Omega_2$ with $\partial
\Gamma_2 =
\partial \Sigma_2'$. Since $\partial \Gamma_1 , \partial \Gamma_2$ are
linked in $\Omega_1 , \Omega_2$ with (segments of)
$\gamma_{y_1,y_3} ,\gamma_{y_3,y_2'}$, respectively, we get
components $\Gamma_i'$ of $B_{r_0-4h} (y) \cap \Gamma_i$ with
$z^{\Gamma}_1 \in \Gamma_1' \cap \gamma_{y_1 , y_3}$ and
$z^{\Gamma}_2 \in \Gamma_2' \cap \gamma_{y_3 , y_2'}$. By Lemma
\ref{l:lone}, each $\Gamma_i'$ is a graph of a function $u_i$ with
$|\nabla u_i| \leq C_g \, \ba$. Hence, since $1 + C_g^2 \, \ba^2 <
1 / \cos^2 \theta_0$,
\begin{equation} \label{e:nv2}
\Gamma_i' \setminus \{ z_i^{\Gamma} \} \subset N_{\theta_0}
(z^{\Gamma}_i ) \, .
\end{equation}
By \eqr{e:nv}, $\gamma_{y_1 , y_2'} \cap N_{\theta_0}
(z^{\Gamma}_i ) = \emptyset$, so \eqr{e:nv2} implies that
$\Gamma_i' \cap \gamma_{y_1 , y_2'} = \{ z^{\Gamma}_i \}$. In
particular, $y_1 , y_2 , y_3$ are in distinct components of
$B_{r_0 - 4 \, h} \setminus (\Gamma_1' \cup \Gamma_2')$. This
completes the second step.

Set $\hat{y} =
\partial B_{r_0 + 10 \, h} \cap \gamma_{0,y}$. Let
$\hat{\gamma}$ be the component of $B_{r_0 + 10 \, h} \cap \gamma$
with $B_{r_0} \cap \hat{\gamma} \ne \emptyset$. Then $\partial
\hat{\gamma} = \{ \hat{y}_1 , \hat{y}_2 \}$ with $\hat{y}_i \in
B_h (\hat{y}) \cap \Sigma_i'$.

The third step is to solve the Plateau problem with $\gamma_3$
together with part of $\partial \Sigma \subset \partial B_{4 \,
r_0}$ as the boundary to get a stable disk $\Gamma_3
\subset B_{4 r_0} \setminus \Sigma$ passing between $\hat{y}_1 ,
\hat{y}_2$. To do this, note that the curve $\gamma_3$ divides the
disk $\Sigma$ into two sub-disks $\Sigma_3^+ , \Sigma_3^-$. Let
$\Omega^+, \Omega^-$ be the components of $B_{4\, r_0} \setminus
(\Sigma \cup \Gamma_1 \cup \Gamma_2)$ with $\gamma_3 \subset
\partial \Omega^+ \cap \partial \Omega^-$. Note that $\Omega^+,
\Omega^-$ are mean convex in the sense of \cite{MeYa} since
$\partial \Gamma_1 \cup \partial \Gamma_2 \subset \Sigma$ and
$\partial \Sigma \subset \partial B_{4 \, r_0}$. Using the first
step, we can label $\Omega^+, \Omega^-$ so $z^+ , z^-$ do not
connect in $H \cap \Omega^+$. By \cite{MeYa}, we get a stable
embedded disk $\Gamma_3 \subset \Omega^+$ with $\partial
\Gamma_3 = \partial \Sigma_3^+$. Using the almost
monotonicity, $\partial B_{r_0 - 10 \, h}(y) \cap \partial
\Gamma_3$ consists of an odd number of points in each of
$\Sigma_{H}^+ , \, \Sigma_{H}^-$. Hence, there is a curve
$\gamma_{+}^{-} \subset \partial B_{r_0 - 10 \, h}(y) \cap
\Gamma_3$ from $\Sigma_{H}^+$ to $\Sigma_{H}^-$. By construction,
$\gamma_{+}^{-} \setminus B_{8 \, h}(y_0) \ne \emptyset$. Hence,
since $\partial B_{r_0 - 10 \, h}(y) \cap T_h (\partial \Gamma_3)
\subset B_{3\, h}(y_0)$, Lemma \ref{l:lone} gives
$\hat{z} \in B_h (\hat{y}_1) \cap \gamma_{+}^{-}$. By the
second step, $\Gamma_3$ is between $\Gamma_1'$ and
$\Gamma_2'$.

Let $\hat{\Gamma}_3$ be the component of $B_{r_0 + 19 \, h} \cap
\Gamma_3$ with $\hat{z} \in \hat{\Gamma}_3$. By Lemma
\ref{l:lone}, $\hat{\Gamma}_3$ is a graph. Finally, since
$\hat{\gamma} \subset B_{r_0 + 10 \, h}$ and $\hat{\Gamma}_3$
passes between $\partial \hat{\gamma}$, this forces
$\hat{\Gamma}_3$ to intersect $\hat{\gamma}$. This contradiction
completes the proof.
\end{proof}

\part{Estimates for stable annuli with slits} \label{p:p2}

In this part, we will show that certain stable disks
starting off as multi-valued graphs remain the same (see Theorem
\ref{t:spin4ever}  below). This is needed in Part \ref{p:p3}
when we generalize the results of Part \ref{p:p1} to
when the surface is not anymore in a slab and in
Part \ref{p:p4} when we show Theorem \ref{t:spin4ever2}.

\begin{Thm} \label{t:spin4ever}
Given $\tau  > 0$, there exist $N_1 , \Omega_1, \epsilon > 0$ so:
Let $\Omega_1 \, r_0 < 1 < R_0 / \Omega_1$, $\Sigma \subset
B_{R_0}$ be a stable embedded minimal disk, $\partial \Sigma
\subset B_{r_0} \cup
\partial B_{R_0} \cup \{ x_1 = 0 \}$, $\partial  \Sigma \setminus \partial B_{R_0}$
is connected.  If $\Sigma$ contains an $N_1$-valued graph
$\Sigma_g$ over $D_1 \setminus D_{r_0}$ with gradient $\leq
\epsilon$, $\Pi^{-1}( D_{r_0}) \cap \Sigma^M \subset \{ |x_3| \leq
\epsilon \, r_0 \}$, and a curve $\eta \subset \Pi^{-1} (D_{r_0})
\cap \Sigma \setminus \partial B_{R_0}$ connects $\Sigma_g$ to
$\partial \Sigma \setminus
\partial B_{R_0}$, then $\Sigma$ contains a $2$-valued graph
$\Sigma_d$ over $D_{R_0/\Omega_1} \setminus D_{ r_0}$ with
gradient $\leq \tau$.
\end{Thm}

Two analytical results go into the proof of this extension
theorem. First, we show that if an almost flat multi-valued graph
sits inside a stable disk, then the outward defined intrinsic
sector from a curve which is a multi-valued graph over a circle
has a subsector which is almost flat (see Corollary \ref{c:stable}
below).  As the  initial multi-valued graph becomes flatter and
the number of sheets in it go up, the subsector becomes flatter.
The second analytical result that we will need is that in a
multi-valued minimal graph the distance between the sheets grows
sublinearly (Proposition \ref{l:grades1}).

After establishing these two facts, the first application
(Corollary \ref{c:cep}) is to extending the middle sheet as a
multi-valued graph.  This is done by dividing the initial
multi-valued graph (or curve in the graph that is itself a
multi-valued graph over the circle) into three parts where the
middle sheet is the second part.  The idea is then that the first
and third parts have subsectors which are almost flat multi-valued
graphs and the middle part (which has curvature estimates since it
is stable) is sandwiched between the two others.  Hence its sector
is also almost flat.

A thing that adds some technical complications to the above
picture is that in the analytical result about almost flat
subsectors it is important that the ratio between the size of the
initial multi-valued graph and how far one can go out is fixed.
This is because the estimate for the subsector comes from a total
curvature estimate which is in terms of this ratio (see
\eqr{e:cbound}) and can only be made small by looking at a fixed
large number of rotations for the graph. This forces us to
successively extend the multi-valued graph.  The issue is then to
make sure that as we move out in the sector and repeat the
argument we have essentially not lost sheets.  This is taken care
of by using the sublinear growth of the separation between the
sheets together with the Harnack inequality (Lemma \ref{l:harnie})
and the maximum principle (Corollary \ref{c:cep}). (The maximum
principle is used to make sure that as we try to recover sheets
after we have moved out then we don't hit the boundary of the disk
before we have recovered essentially all of the sheets that we
started with.) The last thing is a result from \cite{CM3}
to guarantee that as we patch together these multi-valued
graphs coming from different scales then the surface that we get
is still a multi-valued graph over a fixed plane.

Unless otherwise stated in this part, $\Sigma$ will be a stable embedded
disk. Let $\gamma \subset \Sigma$ be a simple curve with unit
normal $\nn_{\gamma}$ and geodesic curvature $k_g$ (with respect
to $\nn_{\gamma}$). We will always assume that $\gamma'$ does not
vanish.
Given $R_1 > 0$, we define the intrinsic sector, see fig. \ref{f:6},
\begin{equation}		\label{e:figf6}
S_{R_1}(\gamma) = \cup_{x \in \gamma} \gamma_x \, ,
\end{equation}
where $\gamma_x$ is the (intrinsic) geodesic starting at $x \in
\gamma$, of length $R_1$, and initial direction $\nn_{\gamma}(x)$.
For $0 < r_1 < R_1$, set $S_{r_1, R_1}(\gamma) = S_{R_1}(\gamma)
\setminus S_{ r_1}(\gamma)$ and $\rho (x) =
\dist_{S_{R_1}(\gamma)} ( x , \gamma)$. For example, if $\gamma =
\partial D_{r_1} \subset \RR^2$ and $\nn_{\gamma} (x) = x / |x|$,
then $S_{r_2,R_1}$ is the annulus
$D_{R_1 +r_1} \setminus D_{r_2+r_1}$.

\begin{figure}[htbp]
    \setlength{\captionindent}{20pt}
    \begin{minipage}[t]{0.5\textwidth}
    \centering\input{shn6.pstex_t}
    \caption{An intrinsic sector over a curve $\gamma$ defined in 
	\eqr{e:figf6}.}
    \label{f:6}    \end{minipage}\begin{minipage}[t]{0.5\textwidth}
\centering\input{shn6a.pstex_t}
    \caption{The curve $\gampa$ containing $\gamma$ goes to $\partial \Sigma$.
($\gampa \setminus \gamma$ is dotted.)}
\label{f:7}    \end{minipage}
\end{figure}

Note that if $k_g > 0$, $S_{R_1} (\gamma) \cap
\partial \Sigma = \emptyset$, and there is a simple curve
$\gampa \subset \Sigma$ with $\gamma \subset \gampa$, $\partial
\gampa \subset \partial \Sigma$, and $\gamma_x \cap \gampa = \{ x
\}$ for any $\gamma_x$ as above (see fig. \ref{f:7}), then the normal exponential map
from $\gamma$ (in direction $\nn_{\gamma}$) gives a diffeomorphism
to $S_{R_1} (\gamma)$. Namely, by the Gauss-Bonnet theorem, an
$n$-gon in $\Sigma$ with concave sides and $n$ interior angles
$\alpha_i > 0$ has
\begin{equation} \label{e:gba}
(n-2) \, \pi \geq \sum_{i=1}^n \alpha_i - \int k_g \geq \sum_{i=1}^n
\alpha_i \, .
\end{equation}
In particular, $n >2$ always and if $\sum_i \alpha_i > \pi$, then
$n > 3$. Fix $x,y \in \gamma$ and geodesics $\gamma_x , \gamma_y$
as above. If $\gamma_x$ had a self-intersection, then it would
contain a simple geodesic loop, contradicting \eqr{e:gba}.
Similarly, if $\gamma_x$ were to intersect $\gamma_y$, then we
would get a concave triangle with $\alpha_1 = \alpha_2 = \pi / 2$
(since $\gamma_x , \gamma_y$ don't cross $\gampa$), contradicting
\eqr{e:gba}.

Note also that $S_{r_1,R_1}(\gamma) = S_{ R_1 -r_1}
(S_{r_1,r_1}(\gamma))$ for $0<r_1<R_1$.

\section{Almost flat subsectors}

We will next show that certain stable sectors contain almost flat
subsectors.

\begin{Lem} \label{l:stable}
Let $\gamma \subset \Sigma$ be a curve with
${\text{Length}}(\gamma) \leq 3\, \pi \, m \, r_1$, $0 < k_g < 2 /
r_1$, $\dist_{\Sigma} (S_{R_1} (\gamma) ,
\partial \Sigma ) \geq r_1 / 2$, $R_1 > 2 \, r_1$.  If there is a simple curve $\gampa
\subset \Sigma$ with $\gamma \subset \gampa$, $\partial \gampa
\subset \partial \Sigma$, and $\gamma_x \cap \gampa = \{ x \}$ for
$x \in \gamma$, then for  $\Omega > 2$ and $2 \, r_1 \leq t \leq 3 R_1 / 4$
\begin{align}
\int_{S_{\Omega r_1 , R_1 / \Omega}(\gamma)} |A|^2 &\leq C_1 \, R_1
/ r_1 + C_2 \, m / \log \Omega \, ,\label{e:cbound}\\
t \, \int_{\gamma} k_g \leq {\text{Length}} (\{ \rho = t \}) &\leq C_3
\, (m +R_1/r_1) \, t \, . \label{e:lbound}
\end{align}
\end{Lem}

\begin{proof}
The boundary of $S_{ R_1} = S_{ R_1} (\gamma)$ has four pieces:
$\gamma$, $\{ \rho = R_1 \}$, and the sides $\gamma_a , \gamma_b$.
Set
\begin{align} \label{e:ellt}
\ell(t) &= {\text{Length}} \, ( \{ \rho = t \} ) \, ,\\ K(t) & =
\int_{ S_{t} } |A|^2 \, .
\end{align}
Since the exponential map is an embedding, an easy calculation
gives
\begin{equation} \label{e:difft1}
\ell'(t) = \int_{ \{ \rho = t \} } k_g > 0 \, .
\end{equation}
Let $d\mu$ be $1$-dimensional Hausdorff measure on the level sets
of $\rho$. The Jacobi equation gives
\begin{equation} \label{e:difft2a}
\frac{d}{dt} (k_g \, d \mu) = |A|^2 / 2\, d \mu \, .
\end{equation}
Set $\bar{K}(t) = \int_{0}^{t} K(s) \, ds$. Integrating
\eqr{e:difft2a} twice, \eqr{e:difft1} yields
\begin{equation}
 \ell(t) = \ell (0) + \int_0^t \left( \int_{\gamma} k_g + K(s) / 2 \right) \,
 ds =
{\text{Length}} ( \gamma) + t \, \int_{\gamma} k_g + \bar{K}(t) /
2 \, . \label{e:difft2bbb}
\end{equation}
This gives the first inequality in \eqr{e:lbound}. Again by the
coarea formula, \eqr{e:difft2bbb} gives
\begin{align} \label{e:difft2b}
R_1^{-2} \, \Area ( S_{R_1} ) &= R_1^{-2} \int_{0}^{R_1} \ell(t)
\leq R_1^{-1} \, {\text{Length}}(\gamma) + \int_{\gamma} k_g / 2 +
R_1^{-2} \int_{0 }^{R_1} \bar{K}(t) / 2 \notag \\ & \leq 6 \, \pi
\, m + R_1^{-2} \int_{0}^{R_1} \bar{K}(t) / 2 \, ,
\end{align}
where the last inequality used $k_g < 2/ r_1$ on $\gamma$,
${\text{Length}}(\gamma) \leq 3 \, \pi \, m \, r_1$, and $R_1 > 2
\, r_1$.

Define a function $\psi$ on $S_{R_1 }$ by $\psi =\psi (\rho ) =
1- \rho /R_1$ and set $d_S= \dist_{\Sigma} ( \cdot , \gamma_a
\cup \gamma_b)$. Define functions $\chi_1 , \chi_2$ on $S_{R_1 }$
by
\begin{align}
\chi_1 =& \chi_1 (d_S ) =
\begin{cases}
d_S /r_1 &\hbox{ if } 0 \leq d_S \leq r_1 \, , \\ 1 &\hbox{
otherwise} \, ,\\
\end{cases}
\\
\chi_2 =&\chi_2 (\rho ) =
\begin{cases}
\rho / r_1 &\hbox{ if } 0 \leq \rho \leq r_1 \, , \\ 1 &\hbox{
otherwise} \, .\\
\end{cases}
\end{align}
Set $\chi = \chi_1 \, \chi_2$. Using $|A|^2 \leq C \, r_1^{-2}$
(by \cite{Sc}) and standard comparison theorems to bound the area
of a tubular neighborhood of the boundary, we get
\begin{align}
\Area (S_{R_1} \cap \{ \chi < 1 \} ) & \leq \cchi\, (R_1 \, r_1 +
m \, r_1^2) \, , \label{e:chiB1} \\ E (\chi_1) + \int_{S_{R_1}
\cap \{ \chi_1 < 1 \}} |A|^2 & \leq \cchi\, R_1 / r_1 \, ,
\label{e:chiB3a} \\ E (\chi) + \int_{S_{R_1} \cap \{ \chi < 1 \}}
|A|^2 & \leq \cchi\, ( R_1 / r_1 + m ) \, . \label{e:chiB3}
\end{align}

Substituting $\chi \psi$ into the stability inequality, the
Cauchy-Schwarz inequality and \eqr{e:chiB3} give
\begin{align} \label{e:stab1}
\int |A|^2 \chi^2 \psi^2 & \leq \int |\nabla ( \chi \psi)|^2 =
\int \, \left( \chi^2 |\nabla \psi|^2 + 2 \chi \, \psi \langle
\nabla \chi , \nabla \psi \rangle + \psi^2 |\nabla \chi|^2 \right)
\notag \\ & \leq 2 \int \chi^2 |\nabla \psi|^2 + 2 \, \cchi(R_1 /
r_1 + m ) \, .
\end{align}
Using \eqr{e:chiB3} and the coarea formula, we have
\begin{equation} \label{e:dotheA}
\int_{0}^{R_1} \psi^2(t) \, K'(t) = \int_{S_{R_1}} |A|^2 \psi^2
\leq \int |A|^2 \chi^2 \psi^2 + \cchi\,( R_1 / r_1 + m ) \, .
\end{equation}
Integrating by parts twice in \eqr{e:dotheA}, \eqr{e:stab1} gives
\begin{align} \label{e:stab3}
2 \, R_1^{-2} \int_{0 }^{R_1} \bar{K}(t) &= \int_{0 }^{R_1}
\bar{K}(t) (\psi^2)'' =
-
\int_{0 }^{R_1} K(t) (\psi^2)' \\ & = \int_{0 }^{R_1} \psi^2 K'(t)
\leq 3 \, \cchi\,( R_1 /r_1 + m) + 2 R_1^{-2} \int_{0 }^{R_1}
\ell(t) \, . \notag
\end{align}
Note that all integrals in \eqr{e:stab3} are in one variable and
there is a slight abuse of notation in regarding $\psi$ as a
function on both $[0,R_1]$ and $S_{R_1}$. Substituting
\eqr{e:difft2b}, \eqr{e:stab3} gives
\begin{equation} \label{e:stab6}
4 \, R_1^{-2} \int_{0 }^{R_1} \ell(t) \leq 24 \, \pi \, m + 3 \,
\cchi\,( R_1 /r_1 + m) + 2 \, R_1^{-2} \int_{0 }^{R_1} \ell(t) \,
.
\end{equation}
In particular, \eqr{e:stab6} gives
\begin{equation} \label{e:stab7}
R_1^{-2} \Area ( S_{R_1}) \leq C_4 \, ( R_1 / r_1 + m) \, .
\end{equation}
Since $\ell (t)$ is monotone increasing (by \eqr{e:difft1}),
\eqr{e:stab7} gives the second inequality in \eqr{e:lbound} for
$t= 3 \, R_1 / 4$. Since the above argument applies with $R_1$
replaced by $t$ where $2 \, r_1 < t < R_1$, we get \eqr{e:lbound}
for $2 \, r_1 \leq t\leq 3 \, R_1 / 4$.

To complete the proof, we will use the stability inequality
together with the logarithmic cutoff trick to take advantage of
the quadratic area growth. Define a cutoff function $\psi_1$ by
\begin{equation}
\psi_1 =\psi_1 (\rho ) =
\begin{cases}
\log ( \rho / r_1) / \log \Omega &\hbox{ on } S_{r_1 , \Omega \,
r_1} \, , \\ 1 &\hbox{ on } S_{\Omega \, r_1 , R_1 / \Omega} \,
,\\ - \log ( \rho / R_1) / \log \Omega &\hbox{ on } S_{R_1 /
\Omega ,R_1} \, ,\\ 0 &\hbox{ otherwise} \, .\\
\end{cases}
\end{equation}
Using \eqr{e:lbound} and \eqr{e:stab7}, we get
\begin{equation}
E(\psi_1) \leq C (m + R_1 / r_1) / \log \Omega \, .
\end{equation}
As in \eqr{e:stab1}, we apply the stability inequality to $\chi_1
\psi_1$ to get
\begin{equation} \label{e:stab1too}
\int |A|^2 \chi_1^2 \psi_1^2 \leq 2 E (\psi_1) + 2 E (\chi_1) \leq
2 \, C (m + R_1 / r_1) / \log \Omega + 2 \cchi\, R_1 / r_1 \, .
\end{equation}
Combining \eqr{e:chiB3a} and \eqr{e:stab1too} completes the proof.
\end{proof}

Using Lemma \ref{l:stable}, we show that large stable sectors have almost flat
subsectors:

\begin{Cor} \label{c:stable}
Given $\omega > 8, 1 > \epsilon > 0$, there exist $m_1 ,
\Omega_1$ so: Suppose $\gamma \subset B_{2 \, r_1} \cap \Sigma$ is a
curve with $1 / (2
\, r_1) < k_g < 2 / r_1$, ${\text{Length}}(\gamma) = 32 \, \pi \,
m_1  \, r_1$, $\dist_{\Sigma} (S_{\Omega_1^2 \, \omega \, r_1 } (\gamma) ,
\partial \Sigma ) \geq r_1 / 2$.  If
there is a simple curve $\gampa \subset \Sigma$ with $\gamma
\subset \gampa$, $\partial \gampa \subset \partial \Sigma$, and
$\gamma_x \cap \gampa = \{ x \}$ for $x \in \gamma$, then
(after rotating $\RR^3$) $S_{ \Omega_1^2 \, \omega \, r_1}
(\gamma)$ contains a $2$-valued graph $\Sigma_d$ over $D_{2\, \omega \,
\Omega_1 \, r_1} \setminus D_{\Omega_1 \, r_1 / 2}$ with gradient
$\leq \epsilon/2$, $|A| \leq \epsilon / (2 \, r)$, and
$\dist_{S_{ \Omega_1^2 \, \omega \, r_1}
(\gamma)} ( \gamma , \Sigma_d ) < 2 \, \Omega_1 \, r_1$.
\end{Cor}

\begin{proof}
We will choose $\Omega_1 > 12$ and then set $m_1 = \omega \,
\Omega_1^2 \, \log \Omega_1$. By Lemma \ref{l:stable} (with
$\Omega = \Omega_1 / 6$, $R_1 = \Omega_1^2 \, \omega \, r_1$, and
$m = 32 \, m_1 / 3$),
\begin{equation} \label{e:fromcb}
\int_{S_{\Omega_1 \, r_1 / 6 , 6\, \Omega_1 \, \omega \,
r_1}(\gamma)} |A|^2 \leq C( \Omega_1^2 \, \omega + m_1 / \log
\Omega_1 ) = 2 \, C \, m_1 / \log \Omega_1 \, .
\end{equation}
Fix $m_1$ disjoint  curves $\gamma_1 , \dots ,
\gamma_{m_1} \subset \gamma$ with ${\text{Length}}(\gamma_i) = 32
\, \pi \, r_1$.  By \eqr{e:fromcb} and since the $S_{
\Omega_1^2 \, \omega \, r_1}(\gamma_i)$ are pairwise disjoint,
there exists $\gamma_i$ with
\begin{equation} \label{e:fromcb2}
\int_{S_{\Omega_1 \, r_1 / 6 , 6\, \Omega_1 \, \omega \,
r_1}(\gamma_i)} |A|^2 \leq 2 \, C / \log \Omega_1 \, .
\end{equation}

To deduce the corollary from \eqr{e:fromcb2} we need a few
standard facts. First, define a map $\Phi : [0,\Omega_1^2\omega
r_1] \times_{\rho /(2 \, r_1) +1}
[0,\text{Length}(\gamma)] \to \Sigma$ by $\Phi
(\rho,x)=\gamma_x(\rho)$. By the Riccati comparison argument
(using $\K_{\Sigma} \leq 0$ and $k_g
> 1 / (2 \, r_1)$ on $\gamma$),
\begin{equation} \label{e:stcoth}
\Phi \text{ is distance nondecreasing and } k_g > \frac{1}{\rho +
2\, r_1} \, .
\end{equation}

Second,
let $\gamma_i / 2 \subset \gamma_i$ be the  subcurve of
length $16 \, \pi \, r_1$ with
$\dist_{\gamma}(\gamma_i / 2 , \partial \gamma_i)
= 8 \, \pi \, r_1$.
Since $k_g > 1 / (2 \, r_1)$ on $\gamma$, we have
$\int_{\gamma_i/2} k_g > 8 \, \pi$. By \eqr{e:difft2a},
$\int_{S_{\Omega_1^2 \, \omega \, r_1} (\gamma_i/2) \cap \{ \rho = t
\}} k_g$ is monotone nondecreasing. In particular, we can choose a
 curve $\tilde{\gamma} \subset \gamma_i/2$ with
\begin{equation} \label{e:deftg}
\int_{S_{ \Omega_1^2 \, \omega \, r_1} (\tilde{\gamma}) \cap \{
\rho = \Omega_1 \, r_1 / 3 \} } k_g = 8 \, \pi \, .
\end{equation}
Set $S = S_{\Omega_1 \, r_1 / 3 , 3\, \Omega_1 \, \omega \,
r_1}(\tilde{\gamma})$ and $\hat{\gamma} = S \cap \{ \rho =
\Omega_1 \, r_1 / 3 \}$.

Third, by the Gauss-Bonnet theorem, \eqr{e:fromcb2}, and \eqr{e:deftg},
(for  $\Omega_1$ large)
\begin{equation} \label{e:sectime}
8 \, \pi \leq \int_{S \cap \{ \rho = t \} } k_g \leq 8 \, \pi +
\int_S |A|^2 / 2 \leq 8 \, \pi + C / \log \Omega_1 \leq 9 \, \pi
\, .
\end{equation}
Note also that, by \eqr{e:stcoth} and \eqr{e:sectime},
${\text{Length}}( S \cap \{ \rho = t \} ) \leq 9 \, \pi \, (t + 2
\, r_1) \leq 14 \, \pi \, t$.

Finally, observe that, by stability,
\eqr{e:fromcb2}, and using \eqr{e:stcoth}, the meanvalue theorem
gives for $y\in S$
\begin{equation} \label{e:pcnd}
\sup_{\cB_{\rho(y)/3}(y)} |A|^2 \leq C_1 \, \rho^{-2} (y)/ \log
\Omega_1 \, .
\end{equation}
Integrating \eqr{e:pcnd} along rays and level sets of $\rho$, we
get
\begin{equation} \label{e:osbound}
\max_{x,y \in S} \dist_{\SS^2} ( \nn(x) , \nn(y) ) \leq C_2 \, (
\log \omega + 1) / \sqrt{\log \Omega_1} \, .
\end{equation}

We can now combine these facts to get the
corollary.
Choose $\Omega_1$ so that $C_2 \, ( \log \omega + 1) / \sqrt{\log
\Omega_1} < C_3 \, \epsilon$. For $C_3$ small, after rotating
$\RR^3$, $S$ is locally a graph over $\{ x_3 = 0\}$ with gradient
$\leq \epsilon/2$. Since $\tilde{\gamma} \subset B_{2 \,
r_1}$ and $\Omega_1 > 12$, we have $\hat{\gamma}
\subset B_{2 \, r_1 + \Omega_1 \, r_1 / 3}
\subset B_{\Omega_1 \, r_1 / 2}$.
Choosing $\Omega_1$ even larger and combining \eqr{e:stcoth},
\eqr{e:sectime}, \eqr{e:pcnd}, and \eqr{e:osbound}, we see that
(the orthogonal projection) $\Pi (\hat{\gamma})$ is a convex
planar curve with total curvature at least $7 \, \pi$, so that its
Gauss map covers $\SS^1$ three times. Given $x \in
\tilde{\gamma}$, set $\tilde{\gamma}_x = S \cap \gamma_x$. By
\eqr{e:pcnd}, $\tilde{\gamma}_x$ has total (extrinsic geodesic)
curvature at most $C_2 \, \log \omega / \sqrt{\log \Omega_1} < C_3
\, \epsilon$ and hence $\tilde{\gamma}_x$ lies in a narrow
cone centered on its tangent ray at $\tilde{x} = \tilde{\gamma}_x
\cap \hat{\gamma}$. For $C_3$ small, this implies that
$\tilde{\gamma}_x$ does not rotate and
\begin{equation}
    | \Pi ( \tilde{x}) - \Pi (
    \tilde{\gamma}_x \cap \{ \rho = t\} ) | \geq  9 \, ( t -
    \Omega_1 \, r_1 / 3 ) / 10 \, .
\end{equation}
Hence, $\Pi (\partial
\tilde{\gamma}_x \setminus \{ \tilde{x} \} ) \notin D_{2\, \omega
\, \Omega_1 \, r_1}$ which
gives $\Sigma_d$ and also
$\dist_{S_{ \Omega_1^2 \, \omega \, r_1}
(\gamma)} ( \gamma , \Sigma_d ) < 2\, \Omega_1 \, r_1$.
\end{proof}

\begin{Rem}     \label{r:uselater}
  For convenience, we assumed that $k_g < 2 / r_1$ in Corollary \ref{c:stable}.
This was used only to apply Lemma \ref{l:stable} and it was used
there only to bound $\int_{\gamma} k_g$ in \eqr{e:difft2b}.
\end{Rem}

 Recall that a domain $\Omega$ is $1/2$-stable if and only if, for
all $\phi \in C_0^{0,1}(\Omega)$, we have the $1/2$-stability
inequality:
\begin{equation} \label{e:stabineq12}
1/2 \int |A|^2 \phi^2 \leq  \int |\nabla \phi |^2 \, .
\end{equation}
Note that the interior curvature estimate of \cite{Sc} extends to
$1/2$-stable surfaces.

In light of Remark \ref{r:uselater}, it is easy to get the
following analog of Corollary \ref{c:stable}:

\begin{Cor} \label{c:uselater}
Given $\omega > 8, 1 > \epsilon > 0, C_0$, $N$, there exist $m_1 ,
\Omega_1$ so: Suppose $\Sigma$ is an embedded minimal disk,
$\gamma \subset \partial \cB_{r_1}(y) \subset \Sigma$ is a
curve, $\int_{\gamma} k_g  < C_0 \, m_1$,   ${\text{Length}}(\gamma) =
m_1  \, r_1$.  If $\cT_{r_1 / 8} (S_{\Omega_1^2 \, \omega \, r_1 } (\gamma))$
is $1/2$-stable,    then
(after rotating $\RR^3$) $S_{ \Omega_1^2 \, \omega \, r_1}
(\gamma)$ contains an $N$-valued graph $\Sigma_N$ over $D_{\omega \,
\Omega_1 \, r_1} \setminus D_{\Omega_1 \, r_1}$ with gradient
$\leq \epsilon$, $|A| \leq \epsilon /  r$, and
$\dist_{S_{ \Omega_1^2 \, \omega \, r_1}
(\gamma)} ( \gamma , \Sigma_N ) < 4 \, \Omega_1 \, r_1$.
\end{Cor}

Note that,  in Corollary \ref{c:uselater},
  $k_g \geq 1/r_1$ and the injectivity of the exponential
map both follow immediately from comparison theorems.

\section{The sublinear growth}

\begin{figure}[htbp]
    \setlength{\captionindent}{20pt}
    \begin{minipage}[t]{0.5\textwidth}
    \centering\input{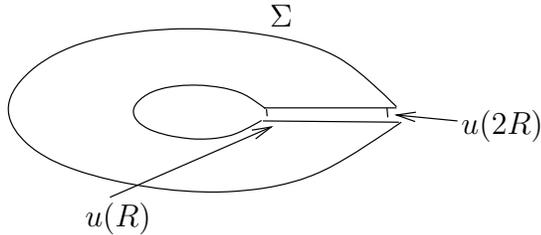}
    \caption{The sublinear growth of the separation $u$ of
the multi-valued graph $\Sigma$:
    $u(2R) \leq 2^{\alpha} \, u(R) $ with $\alpha < 1$.}
    \label{f:8}    \end{minipage}
\end{figure}

This section gives an elementary gradient estimate for
multi-valued minimal graphs which is applied to show that the
separation between the sheets of certain minimal graphs grows
sublinearly; see fig. \ref{f:8}. The example to keep in mind
is the portion of a (rescaled) helicoid in a slab between
two cylinders about the vertical axis. This gives (two)
multi-valued graphs over an annulus; removing a vertical
half-plane through the axis cuts these into sheets
which remain a bounded distance apart.

The next lemma and corollary construct the cutoff function needed
in our gradient estimate.

\begin{Lem} \label{l:graa1}
Given $N > 36/ (1 - \e^{-1/3})^2$, there exists a function $0 \leq \phi \leq 1$ on
$\cP$ with $E(\phi ) \leq 4 \, \pi / \log N$,
\begin{equation} \label{e:phiun2a1}
\phi =
\begin{cases}
1 &\hbox{ if } R/ \e \leq \rho  \leq \e \, R {\text{ and }} |\theta|
\leq 3 \, \pi \, , \\ 0 &\hbox{ if } \rho \leq \e^{-N} \, R
{\text{
or }} \e^{N} \, R \leq \rho
 {\text{
or }} |\theta| \geq \pi \, N \, .\\
\end{cases}
\end{equation}
\end{Lem}

\begin{proof}
After rescaling, we may assume that $R=1$. Since energy is
conformally invariant on surfaces, composing with $z^{3\, N}$
implies that \eqr{e:phiun2a1} is equivalent to $E(\phi ) \leq 4 \, \pi /
\log N$,
\begin{equation} \label{e:phiun2a2}
\phi =
\begin{cases}
1 &\hbox{ if } |\log \rho | < 1 / (3 \, N) {\text{ and }} |\theta| \leq
\pi / N \, ,\\ 0 &\hbox{ if } |\log \rho | > 1/3 {\text{ or }}
|\theta| \geq \pi / 3 \, .\\
\end{cases}
\end{equation}
This is achieved (with $E(\phi) =
2 \, \pi / \log [N (1 - \e^{-1/3})/6]$) by setting
\begin{equation}
\phi =
\begin{cases}
1 &\hbox{ on } \cB_{ 6 / N} (1,0) \, , \\ 1 - \frac{ \log [ N
\, \dist_{\cP} ((1,0),\cdot) / 6 ] }{ \log [ N (1 - \e^{-1/3}) / 6 ] }
&\hbox{ on } \cB_{1- \e^{-1/3} } (1,0) \setminus \cB_{ 6 / N} (1,0) \,
,\\ 0 &\hbox{ otherwise} \, .\\
\end{cases}
\end{equation}
\end{proof}

Given an $N$-valued graph $\Sigma$,
let $\tsma \subset \Sigma$ be the subgraph
(cf.  \eqr{e:defmvg}) over
\begin{equation}    \label{e:subgrn}
    \{ (\rho , \theta) \, | \, r_3 \leq \rho \leq r_4 , \,
        \theta_1 \leq \theta \leq \theta_2 \} \, .
\end{equation}

\begin{Cor} \label{l:graa}
Given $\epsilon_0 , \tau > 0$, there exists $N > 0$ so
if $\Sigma \subset \RR^3$ is an
$N$-valued graph over $D_{\e^{N} \, R} \setminus
D_{\e^{-N} \, R}$ with gradient $\leq \tau$, then there is a cutoff
function $0 \leq \phi \leq 1$ on $\Sigma$ with $E(\phi) \leq
\epsilon_0$, $\phi |_{\partial \Sigma} = 0$, and
\begin{equation}
\phi \equiv 1 {\text{ on }} \tsm  \, . \label{e:phiun2}
\end{equation}
\end{Cor}

\begin{proof}
Since $\tsm \subset \Sigma^{-3\pi,3\pi}_{R/\e , \e R}$
and the projection from $\Sigma$ to $\cP$ is bi-Lipschitz with
bi-Lipschitz constant bounded by $\sqrt{1 + \tau^2}$, the
corollary follows from Lemma \ref{l:graa1}.
\end{proof}

If $u> 0$ is a solution of the Jacobi equation $\Delta u = - |A|^2
u$ on $\Sigma$, then $w= \log u$ satisfies
\begin{equation} \label{e:wlog}
\Delta w = - |\nabla w|^2 - |A|^2 \, .
\end{equation}
The Bochner formula, \eqr{e:wlog}, $\K_{\Sigma} = - |A|^2 / 2$, and
the Cauchy-Schwarz inequality give
\begin{align} \label{e:gwlog}
\Delta |\nabla w|^2 &= 2 \, | \Hess_w |^2 + 2 \langle \nabla w ,
\nabla \Delta w \rangle - |A|^2 \, |\nabla w|^2 \notag \\ &\geq 2
\, | \Hess_w |^2 - 4 \, |\nabla w|^2 \, |\Hess_w| - 4 \, |\nabla
w| \, |A| \, |\nabla A| - |A|^2 \, |\nabla w|^2 \notag \\ &\geq -
2 \, |\nabla w|^4 - 3 \, |A|^2 \, |\nabla w|^2 - 2 \, |\nabla A|^2
\, .
\end{align}
Since the Jacobi equation is the linearization of the minimal
graph equation over $\Sigma$, analogs of \eqr{e:wlog} and
\eqr{e:gwlog} hold for solutions of the minimal graph equation
over $\Sigma$. In particular, standard calculations give
the following analog of \eqr{e:wlog}:

\begin{Lem} \label{l:mge}
There exists $\delta_g > 0$ so if $\Sigma$ is minimal and $u
$ is a positive solution of the minimal graph equation over
$\Sigma$ (i.e., $\{ x + u(x) \, \nn_{\Sigma} (x) \, |\, x \in
\Sigma \}$ is minimal) with $|\nabla u| + |u| \, |A| \leq
\delta_g$, then $w = \log u$ satisfies on $\Sigma$
\begin{equation} \label{e:wlog2}
\Delta w = - |\nabla w|^2 + \dv (a \nabla w) + \langle \nabla w ,
a \nabla w \rangle + \langle b , \nabla w \rangle + (c-1) |A|^2 \,
,
\end{equation}
for functions $a_{ij} , b_j , c$ on $\Sigma$ with $|a| , |c| \leq
3 \, |A| \, |u| + |\nabla u|$ and $|b| \leq 2 \, |A| \, |\nabla
u|$.
\end{Lem}

The following gives an improved gradient estimate, and
consequently an improved bound for the growth of the separation
between the sheets, for multi-valued minimal graphs:

\begin{Pro} \label{l:grades1}
Given $\alpha > 0$, there exist $\delta_p > 0 , N_g > 5 $ so:
Let $\Sigma$ be an $N_g$-valued minimal graph over
$D_{\e^{N_g} \, R} \setminus D_{\e^{-N_g} \, R}$ with gradient
$\leq 1$.  If $0 < u < \delta_p \, R$ is a solution of the minimal graph
equation over $\Sigma$ with $|\nabla u| \leq 1$,
then for $R \leq s \leq 2 \, R$
\begin{align} \label{e:wantit}
\sup_{ \sztp_{R,2R} } |A_{\Sigma}| &+
\sup_{ \sztp_{R,2R} } |\nabla u| / u \leq
\alpha / (4\,R) \, , \\ \label{e:slg}
\sup_{ \sztp_{R,s} } u &\leq (s/R)^{\alpha} \, \sup_{
\sztp_{R,R} } u \, .
\end{align}
\end{Pro}

\begin{proof}
Fix $\epsilon_E > 0$ (to be chosen depending only on $\alpha$).
Corollary \ref{l:graa}
gives $N$ (depending only on $\epsilon_E$) and a  function
$0 \leq \phi \leq 1$ with compact support on
$\Sigma^{-N \pi, N \pi}_{\e^{-N}R , \e^{N} R}$
\begin{equation} \label{e:phiun22}
E(\phi) \leq \epsilon_E {\text{ and }} \phi \equiv 1 {\text{ on }}
\tsm  \, .
\end{equation}
Set $N_g = N+ 1$, so that
$\dist_{\Sigma} ( \Sigma^{-N \pi, N \pi}_{\e^{-N}R , \e^{N} R}
, \partial \Sigma ) > \e^{-N} \, R / 2$ and hence
$|A| \leq C \e^N / R$ on $\Sigma^{-N \pi, N \pi}_{\e^{-N}R , \e^{N} R}$.
Now fix $x \in \sztp_{R,2R}$.  Substituting $\phi$
into the stability inequality,
\eqr{e:phiun22} bounds the total second
fundamental form
of $\tsm$  by
$\epsilon_E$. Hence, by elliptic estimates for the minimal graph
equation,
\begin{equation} \label{e:dlsta2}
\sup_{\cB_{3\, R / 8}(x)} ( R^2 \, |\nabla A_{\Sigma}|^2 +
|A_{\Sigma} |^2 ) \leq C \, \epsilon_E R^{-2} \, .
\end{equation}
Since $\Sigma$ and the graph of $u$ are (locally) graphs with
bounded gradient, it is easy to see that
\begin{equation}    \label{e:smgdp}
    \sup_{\Sigma^{-N \pi, N \pi}_{\e^{-N}R , \e^{N} R} }
    |\nabla u| \leq C \, \e^{N} \, \sup_{\Sigma} |u| / R
    \leq C \, \e^{N} \,  \delta_p \, .
\end{equation}
 Set $w = \log u$. Choose
$\delta_p > 0$ (depending only on $N$), so that \eqr{e:smgdp}
 implies that $w$ satisfies \eqr{e:wlog2} on
$\Sigma^{-N \pi, N \pi}_{\e^{-N}R , \e^{N} R}$
 with $|a| , |b| / |A| , |c| \leq 1 / 4$. Applying
Stokes' theorem to $\dv ( \phi^2 \nabla w - \phi^2 a \nabla w)$
and using the absorbing inequality gives
\begin{equation} \label{e:dlst}
\int_{\cB_{R/2}(x) } |\nabla w|^2 \leq \int \phi^2 |\nabla w|^2
\leq C \, E(\phi) \leq C \, \epsilon_E \, .
\end{equation}
Combining \eqr{e:wlog2} and \eqr{e:dlsta2}, an easy calculation
(as in \eqr{e:gwlog}) shows that on $\cB_{3R/8}(x)$
\begin{equation} \label{e:dlsta4}
\Delta |\nabla w|^2 \geq - C \, |\nabla w|^4 - C \, \epsilon_E \,  R^{-2} \,
|\nabla w|^2 - C \, \epsilon_E \, R^{-4} \, .
\end{equation}
By the rescaling argument of \cite{CiSc} (using the
meanvalue inequality), \eqr{e:dlst} and \eqr{e:dlsta4} imply a
pointwise bound for $|\nabla w|^2$ on $\cB_{R/4}(x)$; combining
this with \eqr{e:dlsta2} gives \eqr{e:wantit} for $\epsilon_E$
 small. Integrating \eqr{e:wantit}  and using that
$(s-R)/R \leq 2 \, \log (s/R)$ gives \eqr{e:slg}.
\end{proof}

\section{Extending multi-valued graphs in stable disks}

Throughout this section $\Sigma \subset B_{R_0}$ is a stable
embedded minimal disk with $\partial \Sigma \subset B_{r_0} \cup
\partial B_{R_0} \cup \{ x_1 = 0 \}$ and
$\partial  \Sigma \setminus \partial B_{R_0}$
  connected.  Fix $0 < \tau_k < 1/4$ so if $\Sigma_g$
is a multi-valued minimal graph over $D_{2R}\setminus D_{R/2}$
with gradient $\leq \tau_k$, then $\Pi^{-1} (\partial D_R) \cap
\Sigma_g$ has geodesic curvature $1 / (2 \, R) < k_g < 2 / R$
(with respect to the outward normal).

The next corollary shows that for
certain such $\Sigma$
containing multi-valued graphs, the middle sheet $\Sigma^M$ extends to
a larger scale.  The main point is to apply Corollary
\ref{c:stable} to get two $2$-valued graphs on a larger scale
with $\Sigma^M$ pinched between them.  We first use the convex hull
property to construct the curves
$\gampa_j$ needed for Corollary
\ref{c:stable}.

\begin{Cor} \label{c:cep}
Given $\omega , m > 1$, $1/4 \geq \epsilon > 0$, there exist
$\Omega_1 , m_0 , \delta$ so  for  $r_0 , r_2 , R_2 , R_0$ with $4
\, \Omega_1 \, r_0 \leq 4 \, \Omega_1 \, r_2 < R_2 < R_0 /(4 \,
\Omega_1 \,\omega)$: Suppose $\Sigma_g\subset \Sigma$ is an $m_0$-valued graph
over $D_{R_2} \setminus D_{r_2}$ with gradient $\leq
\tau_k$, $\Pi^{-1} ( D_{r_2}) \cap \Sigma_g \subset \{ |x_3| \leq
r_2 / 2 \}$,  and separation between the top and bottom sheets
$\leq \delta \, R_2$ over $\partial D_{R_2}$.  If a curve $\eta
\subset \Pi^{-1} (D_{r_2}) \cap \Sigma \setminus
\partial B_{R_0}$ connects $\Sigma_g$ to $\partial \Sigma
\setminus \partial B_{R_0}$, then $\Sigma^M$ extends to an
$m$-valued graph over $D_{\omega \, R_2} \setminus D_{r_2}$ with
gradient $\leq 1$ and $|A| \leq \epsilon / r$ over $D_{\omega \,
R_2} \setminus D_{R_2}$.
\end{Cor}

\begin{proof}
First, we set up the notation. Let $\Omega_1 , m_1 > 1$ be given
by Corollary \ref{c:stable}. Assume that $\Omega_1^2 \, \omega , m
, m_1 \in \ZZ$. Set $m_0 = 24 \, \Omega_1^2 \, \omega + 32 \, m_1
+ m +1$ and $\gamma = \Pi^{-1} (\partial D_{R_2/\Omega_1}) \cap
\Sigma_g$. Since $\Pi^{-1} ( D_{r_2}) \cap \Sigma_g \subset \{
|x_3| \leq r_2 / 2 \}$, the gradient bound gives for $r_2 \leq R
\leq R_2$
\begin{equation} \label{e:controlbd}
\max_{\Pi^{-1} (\partial D_{R}) \cap \Sigma_g} |x_3| \leq r_2 / 2
+ \tau_k \, (R - r_2) \leq R / 2 \, ,
\end{equation}
so that $\gamma \subset B_{2 R_2 / \Omega_1}$. By the definition
of $\tau_k$, $\Omega_1 / (2 \, R_2) < k_g < 2 \, \Omega_1 / R_2$
on $\gamma$.   Arguing on part of
$\Sigma$ itself, by the  convex hull property,
there are $m_0$ components of $\gamma \cap \{
x_1 \geq R_2 / (2 \, \Omega_1) \}$ which are in distinct
components of $\Sigma \cap \{ x_1 \geq R_2 / (2 \, \Omega_1) \}$.
Hence, see fig. \ref{f:9}, 
there are $m_0$ distinct $y_i \in \gamma$ and (nodal)
curves $\sigma_{0} , \dots , \sigma_{m_0-1} \subset \{ x_1 = R_2 /
\Omega_1 \} \cap \Sigma$ with $\partial \sigma_i = \{ y_i , z_i
\}$, $\sigma_i \cap \gamma = \{ y_i \}$, $z_i \in \partial \Sigma
\cap \{ x_1 = R_2 / \Omega_1 \} \subset \partial B_{R_0}$, and for
$i \ne j$
\begin{equation} \label{e:disjtubs}
\dist_{\Sigma} (\sigma_i , \sigma_j) > R_2 / \Omega_1 \, .
\end{equation}
Order the $\sigma_i$'s using the ordering of the $y_i$'s in
$\gamma$ and set $i_1 = 0$, $i_2 = 8 \, \Omega_1^2 \, \omega + 16
\, m_1 $, $i_3 = 16 \, \Omega_1^2 \, \omega + 16 \, m_1  + m$, and
$i_4 = m_0-1$. Let $\gamma_1 , \gamma_2 , \gamma_3 \subset \gamma$
be the curves from $y_{4 \, \Omega_1^2 \, \omega}$ to $y_{4 \,
\Omega_1^2 \, \omega + 16 \, m_1}$, from $y_{12 \, \Omega_1^2 \,
\omega + 16 \, m_1}$ to $y_{12 \, \Omega_1^2 \, \omega + 16 \, m_1
+ m}$, and from $y_{20 \, \Omega_1^2 \, \omega + 16 \, m_1  + m}$
to $y_{20 \, \Omega_1^2 \, \omega + 32 \, m_1  + m}$,
respectively. Hence, $\gamma_1 , \gamma_2 , \gamma_3 \subset
\gamma$ are $16 \, m_1$-, $m$-, $16 \, m_1$-valued graphs,
respectively, with $\gamma_2$ centered on $\Sigma^M$, each
$\gamma_j$ between $y_{i_j}$ and $y_{i_{j+1}}$, and for $j=1,2,3$
\begin{equation} \label{e:lotsis}
 \min_{ \{ k \, | \, y_k \in \gamma_j \} } \,
     \{ | i_j - k| , |i_{j+1} - k|   \} \geq 4 \,
\Omega_1^2 \, \omega \, .
\end{equation}

Next, we  construct the curves $\gampa_j$ needed to apply
Corollary \ref{c:stable}  to  each $\gamma_j$.  We will also use
\eqr{e:disjtubs} and \eqr{e:lotsis} to separate the $\gamma_j$'s.
 For $k_1 < k_2 $, let $\gamma (k_1 , k_2) \subset
\Sigma$ be the union of $\sigma_{k_1}$, $\sigma_{k_{2}}$, and the
curve in $\gamma$ from $y_{k_1}$ to $y_{k_{2}}$. Since $\Sigma$ is
a disk,  $\partial \gamma (k_1 , k_2) \subset
\partial \Sigma$,  and $\partial \Sigma \setminus \partial B_{R_0}$
is connected, one component $\Sigma(k_1,k_2)$ of $\Sigma \setminus
\gamma (k_1 , k_2)$ has $\partial \Sigma(k_1,k_2) \cap \partial
\Sigma \subset \partial B_{R_0}$. Using that the $\sigma_i$'s do
not cross $\eta$, it is easy to see that $\nn_{\gamma}$ points
into $\Sigma (k_1 , k_2) $ and
\begin{equation}    \label{e:inorout}
    \Sigma (j_1 , j_2) \cap \Sigma (k_1 , k_2) =
    \Sigma (\max \{ j_1 , k_1 \} , \min \{ j_2 , k_2 \} ) \, ,
\end{equation}
where, by convention, $\Sigma (k_1 , k_2) = \emptyset$ if $k_1 >
k_2$. Set $\gampa_j = \gamma(i_j , i_{j+1})$ and note that
 $\gamma_j \subset \gampa_j$ and $\partial \gampa_j \subset
\partial \Sigma$.   Set $S_j = S_{ \Omega_1 \,
\omega \, R_2} (\gamma_j)$. By \eqr{e:lotsis} and \eqr{e:inorout},
any curve $\tilde{\eta} \subset \Sigma (i_j , i_{j+1})$
from $\gamma_j$ to $\gampa_j \setminus (\gamma \cup \partial B_{R_0})$
hits at least $4 \, \Omega_1^2 \, \omega$ of the $\sigma_i$'s and
so, by \eqr{e:disjtubs}, ${\text{Length}}(\tilde{\eta}) > 2 \,
\Omega_1 \, \omega \, R_2$.
Combining this with $R_0 > 4 \Omega_1 \, \omega \, R_2$, we get
\begin{equation}        \label{e:newdf}
        \dist_{\Sigma(i_j , i_{j+1})} (\gamma_j , \partial \Sigma(i_j , i_{j+1})
        \setminus \gamma_j ) >
        2 \, \Omega_1 \, \omega \, R_2 \, .
\end{equation}
Fix $x \in \gamma_j$ and $\gamma_x$
(the geodesic normal to $\gamma_j$ at $x$ and of length $\Omega_1
\, \omega \, R_2$). By \eqr{e:gba}, the first point (after $x$)
where $\gamma_x$ hits $\partial \Sigma (i_j , i_{j+1})$ cannot be
in $\gamma$. Consequently, \eqr{e:newdf} implies that
$\gamma_x \subset \Sigma (i_j , i_{j+1})$ so
$\gamma_x \cap \gampa_j = \{ x \}$ and $\gampa_j$ separates $S_j$
from $S_k \cup \cT_{R_2/(2\, \Omega_1)}(\partial \Sigma)$ for
$j\ne k$.

The rest of the proof (see fig. \ref{f:10})
is to sandwich $\Sigma^M$ between two graphs
that will be given by Corollary \ref{c:stable} and then deduce
from stability that $\Sigma^M$ itself extends to a graph. Namely,
applying Corollary \ref{c:stable} to $\gamma_1 , \gamma_3$ (with
$r_1 = R_2 / \Omega_1$), we get $2$-valued graphs $\Sigma_{d,1}
\subset S_1$, $\Sigma_{d,3} \subset S_3$ over $B_{2\, \omega \,
R_2} \cap P_i \setminus B_{R_2/ 2}$ ($i=1,3$) with $|A| \leq
\epsilon / ( 2 \, r)$ and gradient $\leq \epsilon/2 \leq 1/8$.
Here  $P_i$ is a plane through $0$.  Using $|A| \leq \epsilon / ( 2 \, r)$
and $\dist_{S_i}(\gamma , \Sigma_{d,i}) < 2 \, R_2$,
it is easy to see that $\Sigma_{d,i}
\cap \Sigma_g \ne \emptyset$. Hence, $\Sigma_{d,i}$ contains a
$3/2$-valued graph $\Sigma_i$ over $D_{3\, \omega \, R_2/2}
\setminus D_{2 \, R_2/ 3}$ with gradient  $\leq \tan \left(
\tan^{-1} (1/4) + 2 \tan^{-1} (1/8) \right) < 3/4$.   By
construction, $\Sigma^M$ is pinched between $\Sigma_1, \Sigma_3$
which are graphs over each other with separation
$\leq \omega^C \, \delta \, R_2$ (by the Harnack
inequality). Since $\Sigma$ is stable, it follows that if $\delta$
is small, then $\Sigma^M$ extends to an $m$-valued graph
$\Sigma_2$ over $D_{5 \, \omega \, R_2 / 4 } \setminus D_{4 \, R_2
/ 5}$ with $\Sigma_2$ between $\Sigma_1$ and $\Sigma_3$. In
particular, $\Sigma_2$ is a graph over $\Sigma_1$. Finally, using
that $\Sigma_1$ is a graph with gradient  $\leq 3/4$ and $|A| \leq
\epsilon / (2\, r)$, we get that $\Sigma_2$ is a graph with
gradient  $\leq 1$  and $|A| \leq \epsilon / r$ (cf. Lemma
\ref{l:lone}).
\end{proof}

\begin{figure}[htbp]
    \setlength{\captionindent}{20pt}
    \begin{minipage}[t]{0.5\textwidth}
    \centering\input{shn8.pstex_t}
    \caption{The proof of Corollary \ref{c:cep}:
	The nodal curves.}
    \label{f:9}    \end{minipage}\begin{minipage}[t]{0.5\textwidth}
\centering\input{shn9.pstex_t}
    \caption{The proof of Corollary \ref{c:cep}:
	Sandwiching between two graphical pieces.}
\label{f:10}    \end{minipage}
\end{figure}

Combining this and Proposition \ref{l:grades1}, $\Sigma^M$ extends with separation
growing sublinearly:

\begin{Cor} \label{l:ext}
Given $1/4 \geq \epsilon > 0$, there exist $\Omega_0 , m_0 ,
\delta_0 > 0$ so for any $r_0 , r_2 , R_2 , R_0$ with
$\Omega_0 \, r_0 \leq \Omega_0 \, r_2 < R_2 < R_0 / \Omega_0$: Suppose
$\Sigma_g\subset \Sigma$ is an $m_0$-valued graph over $D_{R_2}
\setminus D_{r_2}$ with gradient $\leq \tau_1 \leq \tau_k$,
$\Pi^{-1} (D_{r_2}) \cap \Sigma_g \subset \{ |x_3| \leq r_2 / 2
\}$, and separations between the top and bottom sheets of
$\Sigma^M$($\subset \Sigma_g$) and $\Sigma_g$ are $\leq \delta_1
\, R_2$ and $\leq \delta_0 \, R_2$, respectively, over $\partial
D_{R_2}$.  If a curve $\eta \subset \Pi^{-1} (D_{r_2}) \cap \Sigma
\setminus \partial B_{R_0}$ connects $\Sigma_g$ to $\partial
\Sigma \setminus \partial B_{R_0}$, then $\Sigma^M$ extends as a
graph over $D_{2 \, R_2} \setminus D_{r_2}$ with gradient $\leq
\tau_1 + 3 \, \epsilon$, $|A| \leq \epsilon / r$ over $D_{2 \,
R_2} \setminus D_{R_2}$, and, for $R_2 \leq s \leq 2 \, R_2$,
separation  $\leq (s / R_2)^{1/2} \, \delta_1 \, R_2$ over $D_{s}
\setminus D_{R_2}$.
\end{Cor}

\begin{proof}
Let $\delta_p > 0 , N_g > 5$ be given by Proposition
\ref{l:grades1} with $\alpha = 1/2$.  Let $\Omega_1 , m_0 , \delta
> 0$ be given by Corollary \ref{c:cep} with $m = N_g + 3$ and
$\omega = 2 \, \e^{N_g}$.   We will set $\delta_0 = \delta_0
(\delta , \delta_p , N_g)$ with $\delta > \delta_0 > 0$ and
$\Omega_0 = 4 \, \Omega_1 \, \e^{N_g}$. By Corollary \ref{c:cep},
$\Sigma^M$ extends to a graph
$\Sigma^{-(N_g+3)\pi,(N_g+3)\pi}_{r_2,2 \, \e^{N_g} \, R_2}$ of a
function $v$   with $|\nabla v| \leq 1$ and $|A| \leq \epsilon /
r$ over $D_{2 \, \e^{N_g} \, R_2} \setminus D_{R_2}$. Integrating
$|\nabla |\nabla v|| \leq |A| \, (1 + |\nabla v|^2)^{3/2}\leq
2^{3/2} \, \epsilon / r$, we get
 that $|\nabla v| \leq \tau_1 + 4 \, \epsilon \, \log 2
\leq \tau_1 + 3 \, \epsilon$ on $D_{2R_2} \setminus D_{R_2}$.

    For $\delta_0 = \delta_0 (N_g, \delta_p) > 0$, writing $\Sigma$
as a graph over itself and using the Harnack inequality, we get
 a solution $0 < u < \delta_p \, R_2$ of the minimal graph equation
on an $N_g$-valued graph over $D_{\e^{N_g} \, R_2}
\setminus D_{\e^{-N_g} \, R_2}$.   Applying Proposition
\ref{l:grades1} to $u$ gives the last claim.
\end{proof}

The next lemma uses the Harnack inequality
to show that if
$\Sigma^M$ extends with small separation, then so do the other sheets.
The only complication is to keep track of $\partial \Sigma$.

\begin{Lem} \label{l:harnie}
Given $N \in \ZZ_+$, there exist $C_3 , \delta_2 > 0$ so for
$ r_0 \leq s < R_0/8$: Suppose $\Sigma_g
\subset \Sigma\cap \{ |x_3| \leq 2 \,s \}$ is an $N$-valued graph  over $D_{2 \, s} \setminus D_{s}$.  If
a curve $\eta \subset \Pi^{-1} (D_{s}) \cap \Sigma \setminus
\partial B_{R_0}$ connects $\Sigma_g$ to $\partial \Sigma
\setminus \partial B_{R_0}$,
 and $\Sigma^M$
extends graphically over $D_{4\, s} \setminus D_{s}$ with gradient
$\leq \tau_2 \leq 1$ and separation $\leq \delta_3 \, s \leq
\delta_2 \, s$, then $\Sigma_g$ extends to an $N$-valued graph over $D_{3
\, s} \setminus D_{s}$ with gradient $\leq \tau_2 + C_3 \,
\delta_3$ and separation between the top and bottom sheets $\leq
C_3 \, \delta_3 \, s$.
\end{Lem}

\begin{proof}
Suppose $N$ is odd (the even case is virtually
identical).   Fix $y_{-N} ,
\dots, y_{N} \in \Sigma_g$ with
$y_j$ over $\{\rho = 2 \, s,\, \theta = j \, \pi \}$.
 Let  $\gamma_0 , \gamma_2 \subset \Sigma^M$ be the graphs over
$\{ 2s \leq \rho \leq 3s  , \theta = 0 \}$ and
$\{ 2s \leq \rho \leq 3 s  , \theta = 2 \, \pi \}  $,
respectively, with $\partial
\gamma_0 = \{ y_0 , z_0 \}$ and $\partial
\gamma_2 = \{ y_2 , z_2 \}$.

Arguing as
in the proof of Corollary \ref{c:cep}, there are nodal curves
$\sigma_{-N} , \dots , \sigma_{N} \subset \{ x_1 = - 2 \, s
\} \cap \Sigma$ from $y_j$ (for $j$ odd) to $\partial B_{R_0}$ so: (1)
Any curve in
$\Sigma \setminus \Pi^{-1}(\partial D_{2 \, s})$
from $z_0$ to $\partial \Sigma \setminus \partial B_{R_0}$ hits either
every $\sigma_j$ with $j > 0$ or every  $\sigma_j$ with $j < 0$;
 (2) for $i < j$, $\sigma_i$ and $\sigma_j$ do not
connect in $\Pi^{-1} (D_{4s}) \cap \{ x_1 \leq - 2 \, s \} \cap \Sigma$;
and (3) $\dist (\cup_j \sigma_j , \partial \Sigma \setminus \partial
B_{R_0}) \geq s$.
Note that (2) follows easily from the convex hull property when
$i \ne -N$ or $j \ne N$; the case $i= -N$ and $j=N$ follows since
$\Sigma$ separates $y_{-N} , y_N$ in
$\Pi^{-1} (D_{4s}) \cap \{ x_1 \leq - 2 \, s \}$.

By \cite{Sc} and the Harnack inequality for the minimal graph
equation, there exist $C_4 , \delta_4 > 0$ so if $z_3 , z_4
\in \Sigma \setminus \cT_{s/4}(\partial \Sigma)$, $\Pi(z_3) = \Pi
(z_4)$, and $0 < |z_3 - z_4| \leq \delta_5 \, s \leq \delta_4 \,
s$, then $\cB_{s/8}(z_4)$ is a graph over (a subset of)
$\cB_{s/7}(z_3)$ of a function $u > 0$ with $|\nabla u| \leq \min
\{ 1/2 , C_4 \, \delta_5 \}$.
The lemma now follows easily by repeatedly applying this and using (1)--(3)
to stay away from $\partial \Sigma$ until we have recovered all
$N$ sheets.
\end{proof}

\section{Proof of Theorem \ref{t:spin4ever}}

Let again $\Sigma
\subset B_{R_0}$ be a stable embedded disk with $\partial
\Sigma \subset B_{r_0} \cup
\partial B_{R_0} \cup \{ x_1 = 0 \}$ and $\partial  \Sigma \setminus \partial B_{R_0}$
 connected.  We will use the notation of
    \eqr{e:subgrn}, so that $\Sigma^{0,2\pi}_{r_3, r_4}$ is an annulus
with a slit as defined in \cite{CM3}.
An easy consequence of
theorem $3.36$  of \cite{CM3} is:

\begin{Lem} \label{l:336}
Given $  \tau_0 > 0$,
there exists $0 < \epsilon_1=\epsilon_1 ( \tau_0 ) < 1/24$ so: If $2 r_0 \leq 1 < r_3 \leq R_0 / 2$ and
$\Sigma^{0,2\pi}_{1,r_3} \subset \Sigma$ is
the graph of a function $u$
with
$|\nabla u| \leq 1/12$, $\max_{\Sigma^{0,2\pi}_{1,1} }
    (|u| + |\nabla u|) \leq 2\, \epsilon_1$,
$|A| \leq \epsilon_1 / r$, and for $1 \leq t \leq r_3$ the
separation over $\partial D_t$ is
$\leq 4 \, \pi \, \epsilon_1  \, t^{1/2}$, then
$|\nabla u| \leq \tau_0 $.
\end{Lem}

Lemma \ref{l:336} follows from
theorem $3.36$  of \cite{CM3} and
two facts.  First,  since $\Sigma$ is a
graph over a larger set in $\cP$ (using stability and that
$\partial \Sigma \subset B_{r_0} \cup \partial B_{R_0} \cup
\{ x_1 = 0 \}$), the bound for the separation and
estimates for the minimal graph
equation over $\Sigma$ give a bound for the difference in
the two values of $\nabla u$ along the slit (cf. Proposition
\ref{l:grades1}).
Second, theorem $3.36$ of \cite{CM3} actually applies directly to
$B_{3 r_3/4} \cap \Sigma^{0,2\pi}_{1,r_3} \setminus B_{2}$
to get $|\nabla u| \leq \tau_0 / 2$ on $D_{r_3 /2} \setminus
D_{2}$; integrating $|\nabla |\nabla u|| \leq |A| \, (1 + |\nabla u|^2)^{3/2}
\leq 2 \, \epsilon_1 / r$ then gives $|\nabla u| \leq \tau_0 $
on $D_{r_3} \setminus D_{1}$.

We will prove Theorem \ref{t:spin4ever} by repeatedly applying Corollary \ref{l:ext}
to extend $\Sigma^M$ as a graph,
Lemma \ref{l:336} to get an improved gradient bound,
and then Lemma \ref{l:harnie} to extend additional sheets.

\begin{proof}
(of Theorem \ref{t:spin4ever}).
Set $\tau_0 = \min \{ \tau , \tau_k  , 1/24 \} / 2 $ and
let $0 < \epsilon_1=\epsilon_1 (\tau_0 ) < 1/72$
be given by Lemma \ref{l:336}.
$\Omega_0 , m_0
, \delta_0$ be given by Corollary \ref{l:ext} (depending on $\epsilon_1$)
and $C_3 , \delta_2 > 0$ be from
Lemma \ref{l:harnie} with $N = m_0$.
Set $N_1 = m_0$, $\Omega_1 = 2\, \Omega_0$,
 and choose $\epsilon > 0$ so:
\begin{equation}
    \epsilon < \min \, \{  \frac{\epsilon_1}{ 2} , \,
    \frac{\tau_0 }{4 \, \pi \, 2^{1/2} \,  C_3} , \,
    \frac{\delta_0}{2 \, \pi \,2^{1/2} \,  C_3} , \,
    \frac{\delta_0}{2 \, \pi \, m_0} , \,
     \frac{\delta_2}{4 \, \pi \, 2^{1/2}}
      \} \, ,
\end{equation}
$\Pi^{-1} (D_{r_0}) \cap \Sigma_g \subset \{ |x_3| \leq r_0 /2 \}$,
and $|A| \leq \epsilon_1 / r$ on $\Sigma^M \setminus B_{2 \, r_0}$.
To arrange the last condition, we  use the gradient bound, stability,
 and second derivative estimates for the minimal graph
equation (in terms of the gradient bound).
Note that, using gradient $\leq \epsilon$,
the separation between the top and bottom
sheets of $\Sigma^{0,2\pi}_{r_0,1}$ and $\Sigma^{-m_0 \pi , m_0 \pi}_{r_0,1}$
over $\partial D_t$ are at most $2 \, \pi \, \epsilon \, t$ and
$2 \, \pi \, m_0 \, \epsilon \, t$, respectively.  Note also that
$\Pi^{-1} (D_{3 \,r_0}) \cap \Sigma_g \subset
\{ |x_3| \leq 3 \, \epsilon \, r_0 \}$.

(1) Apply Corollary \ref{l:ext} (with $r_2 = r_0 , R_2 = 1,
\tau_1 = 2 \, \tau_0 $)
to extend
$\Sigma^{0,2\pi}_{r_0,1}$ to a graph
$\Sigma^{0,2\pi}_{r_0 , 2}$
with gradient  $\leq 2 \, \tau_0  + 3 \, \epsilon_1< 1/12$,
$|A| \leq \epsilon_1 / r$
over $D_{2} \setminus D_1$, and,
for $1 \leq t \leq 2$,
\begin{equation}    \label{e:sepless}
    {\text{separation }} \leq  2 \, \pi \epsilon \, t^{1/2} {\text{ over }}
     \partial D_{t}  \, .
\end{equation}

(2) By Lemma \ref{l:336}  (with $r_3 =  2$),
$\Sigma^{0,2\pi}_{1 , 2}$
and hence $\Sigma^{0,2\pi}_{r_0 , 2}$ have
 gradient   $\leq \tau_0 $.

(3) By Lemma \ref{l:harnie} (with $N= m_0 , s=1/2 ,
\tau_2 = \tau_0  , \delta_3 =
4 \, \pi \epsilon \, 2^{1/2}$),
$\Sigma^{0,2\pi}_{r_0 , 3/2}$  is contained in an $m_0$-valued graph
$\Sigma^{-m_0 \pi , m_0 \pi}_{r_0 , 3/2} \subset \Sigma$
over $D_{3/2} \setminus D_{r_0}$ with gradient
$\leq \tau_0  + C_3 \, 4 \, \pi \epsilon \, 2^{1/2} < 2 \, \tau_0 $
and separation $\leq C_3 \, 2 \, \pi \epsilon \, 2^{1/2} < \delta_0$.

Repeat (1)--(3) with:
(1) $R_2 = 3/2$ to extend $\Sigma^{0,2\pi}_{r_0,3/2}$ to
$\Sigma^{0,2\pi}_{r_0 , 3}$ with \eqr{e:sepless} holding for
$1\leq t \leq 3$, (2) $r_3 = 3$ so that $\Sigma^{0,2\pi}_{r_0 , 3}$
has gradient $\leq \tau_0 $, (3) $s= 3/2$ to get
$\Sigma^{-m_0 \pi , m_0 \pi}_{r_0 , 9/2} \subset \Sigma$,
and then again
(1) $R_2 = 9/2$, etc., giving the theorem.
\end{proof}

\part{The general case of Theorem \ref{t:tslab}}  \label{p:p3}

\section{Constructing multi-valued graphs in disks in slabs} \label{s:murs}

Using Part \ref{p:p1}, we show next that an embedded
minimal disk in a slab contains a multi-valued graph if it is not a graph.
We can therefore apply Part \ref{p:p2} to get almost flatness of a
corresponding stable disk past the slab. This is
needed when the minimal surface is not in a thin slab.

\begin{Pro} \label{p:cab}
There exists $\beta > 0$ so: If $\Sigma^2 \subset B_{r_0}
\cap \{ | x_3| \leq \beta \, h \}$ is an embedded minimal disk,
$\partial \Sigma \subset
\partial B_{r_0}$, and a component $\Sigma_1$ of
$B_{10 \, h } \cap \Sigma$ is not a graph,
then $\Sigma$ contains an $N$-valued graph over $D_{r_0 - 2 \, h}
\setminus D_{(60 + 20 \, N) \,h}$.
\end{Pro}

\begin{proof}
The proof has four steps. First we show, by using
Lemma \ref{l:la} twice,
that over a truncated sector in the plane, i.e., over
\begin{equation} \label{e:asin2}
S_{s_1,s_2}(\theta_1 , \theta_2) = \{ (\rho ,
\theta ) \, | \, s_1 \leq \rho \leq s_2 , \, \theta_1
\leq \theta \leq \theta_2 \}
\end{equation}
we have $3$ components of $\Sigma$.
Second, we separate these by stable disks
and order them by height.
Third, we use Proposition \ref{p:pslabn1} to
show that the ``middle'' component is a graph over
a large sector.
Fourth, we repeatedly use the appendix to extend
the top and bottom components around the annulus
and then  Proposition \ref{p:pslabn1} to extend
the middle component as a graph.  This will give
the desired multi-valued graph.

\begin{figure}[htbp]
    \setlength{\captionindent}{20pt}
    \begin{minipage}[t]{0.5\textwidth}
    \centering\input{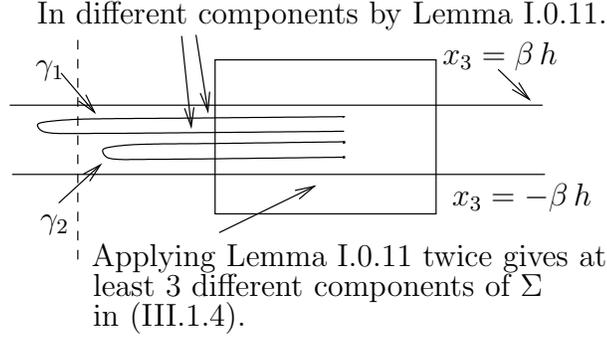}
    \caption{Proof of Proposition \ref{p:cab}:  Step 1:
     Finding the 3 components.}
    \label{f:caba}    \end{minipage}
\end{figure}

  For $j =1,2$, let $\Sigma_j$ be the component of $B_{20 \, j \, h}
\cap \Sigma$ containing $\Sigma_1$. By the maximum principle,
each $\Sigma_j$ is a disk.
Rado's theorem gives
$z_j \in \Pi^{-1} (
\partial D_{(20 \, j - 10) \, h}) \cap \Sigma_j$ for $j= 1 , 2$ where
$\Sigma$ is not graphical (see, e.g., \cite{CM1}).
Rotate $\RR^2$ so
 $z_{1} , z_{2} \in \{ x_1 \geq 0 \}$
and set $z = (r_0 , 0 ,0)$. Apply Lemma \ref{l:la} twice
as in the first step of
the proof of Proposition \ref{p:pslabn1}
to get (see fig. \ref{f:caba}): (1) Disjoint curves 
$\gamma_1 , \gamma_2 \subset
\Sigma$ with $\partial \gamma_k \subset \partial B_{r_0/2}$,
\begin{equation} \label{e:thegammai}
\gamma_k \subset B_{5 \, h }(z_{k}) \cup T_{h }(
\partial D_{(20 \, k - 10) \, h} \cap \{ x_1 \geq 0 \})
\cup T_{ h }(\gamma_{0,z/2}) \, ,
\end{equation}
and which are $C \, \beta \, h$-almost monotone in
$T_{h }(\gamma_{0,z/2})\setminus B_{20 \, k \, h}$.
(2) For $k=1 ,2$
and $y_0 \in \gamma_{0,z/2} \setminus B_{20 \, k \, h}$,
there are components $\Sigma_{y_0,k,1}' \ne
\Sigma_{y_0,k,2}'$ of $B_{5 \, h}(y_0) \cap \Sigma$ each
containing points of $B_{h}(y_0) \cap \gamma_{k}$.
It follows from (2) that, for
$k=1,2$, there are components $\Sigma_{k,1} , \Sigma_{k,2}$ of
$\Pi^{-1} (S_{42 h, r_0 - 2 \, h}(-3\pi /4 , 3 \pi / 4)) \cap \Sigma$
with $\Sigma_{z/2,k,i}' \subset \Sigma_{k,i}$ and which do not connect in
$\Pi^{-1} (S_{40 h , r_0}(-7\pi /8 , 7 \pi / 8)) \cap \Sigma$.
Namely, $\Sigma$ would otherwise contain a disk
violating the maximum principle (as in the second step of Lemma
\ref{l:la}).  The same argument gives
$\Sigma_{i_1,i_1} , \Sigma_{i_2,i_2}, \Sigma_{i_3,i_3}$
which do not connect in
\begin{equation}	\label{e:eqnfor13}
	\Pi^{-1} (S_{40 h , r_0}(-7\pi /8 , 7 \pi / 8)) \cap \Sigma \, .
\end{equation}

By the second step of Proposition
\ref{p:pslabn1}, if $\Sigma_{i,j} ,\Sigma_{k,\ell}$ do not
connect in $\Pi^{-1} (S_{40 h , r_0}(-7\pi /8 , 7 \pi / 8)) \cap \Sigma$,
then there
is a stable embedded disk $\Gamma_{\alpha}$ with $\partial
\Gamma_{\alpha} \subset \Sigma$, $\Gamma_{\alpha} \cap \Sigma =
\emptyset$, and a graph $\Gamma_{\alpha}' \subset \Gamma_{\alpha}$ over
$S_{41 h, r_0 - h}(-13 \pi /16 , 13 \pi / 16)$
separating $\Sigma_{i,j} , \Sigma_{k,\ell}$. Applying this twice
(and reordering the $k_{\ell} , i_{\ell} $),  we get
$\Gamma_{1}' \subset \Gamma_{1}$,
$\Gamma_{2}' \subset \Gamma_{2}$
so each
$\Sigma_{k_{\ell} , i_{\ell} }$ is below $\Gamma_{\ell}'$ which is below
$\Sigma_{k_{\ell+1} , i_{\ell+1}}$.
Let $\gamma^t_1$ and $\gamma^b_1$ be top
and bottom components of
$\cup_j \gamma_{j} \setminus B_{40 \,h }$  intersecting $\partial B_{r_0/2}$.
Since $\Sigma_1 \subset \Sigma_2$, a curve
$\gamma^m_1 \subset B_{40 \, h} \cap \Sigma$
connects $\gamma^t_1$ to $\gamma^b_1$.

\begin{figure}[htbp]
    \setlength{\captionindent}{20pt}
    \begin{minipage}[t]{0.5\textwidth}
\centering\input{shn11b.pstex_t}
    \caption{Proof of Proposition \ref{p:cab}:  Step 3:
     Extending the middle component as a graph.}
\label{f:cabb}    \end{minipage}\begin{minipage}[t]{0.5\textwidth}
    \centering\input{shn11c.pstex_t}
    \caption{Proof of Proposition \ref{p:cab}:  Step 4: Extending the
    top and bottom components by the maximum principle.  They stay disjoint
    since the middle component is a graph separating them.}
    \label{f:cabc}    \end{minipage}
\end{figure}

See fig. \ref{f:cabb}.
By a slight variation of Proposition
\ref{p:pslabn1} (with $\gamma =
\gamma^t_1 \cup \gamma^m_1 \cup \gamma^b_1$),
 the  middle component $\Sigma_{k_2 , i_2}$ is a graph
over
$S_{42 \, h, r_0 - 2 \, h}(-3\pi /4 , 3 \pi / 4)$.
This variation follows from steps one and three of that proof
(step two there constructs barriers $\Gamma_i$
which were constructed here above).

See fig. \ref{f:cabc}.
 Corollary \ref{c:c3} gives curves
$\gamma^t_2 , \gamma^b_2 \subset (B_{44 \, h} \cup T_{h }
(\gamma_{0,(0,r_0 ,0)}) \setminus \Pi^{-1} (D_{42 \, h})) \cap \Sigma$ from
 $\partial B_{43 \, h} \cap \gamma^t_1$ and 
$\partial B_{43 \, h} \cap \gamma^b_1$,
respectively, to $\partial B_{r_0/2}$.  In particular,
 $\gamma^b_2$ is below $\Sigma_{k_2 , i_2}$
and $\gamma^t_2$ is
above $\Sigma_{k_2 , i_2}$; i.e.,
$\Sigma_{k_2 , i_2}$ is still
a middle component.
Again by the maximum  principle, this
gives $3$ distinct components of
$\Pi^{-1}(S_{46 \,h ,r_0 - 2 \, h}(- \pi / 4 , 5 \pi / 4 )) \cap \Sigma$
which do not connect in
$\Pi^{-1}(S_{45 \, h  ,r_0}(-3 \pi / 8 ,11 \pi / 8 )) \cap \Sigma$.
By Proposition
\ref{p:pslabn1}, $\Sigma_{k_2 , i_2}$
further extends as a graph over
$S_{46 \,h ,r_0 - 2 \, h}(- \pi / 4 , 5 \pi / 4 )$, giving a
graph
$\Sigma^{-3\pi /4 , 5\pi /4}_{46 \, h, r_0 - 2 \, h}$
over $S_{46 \,h ,r_0 - 2 \, h}(- 3 \pi / 4 ,5 \pi / 4 )$.
By Rado's theorem,
 this graph cannot close up.  Repeating this with
$\gamma^t_3 , \gamma^b_3 \subset (B_{49 \, h} \cup T_{h }
(\gamma_{0,(-r_0,0 ,0)}) \setminus \Pi^{-1} (D_{47 \, h})) \cap \Sigma$,
etc., eventually gives the proposition.
\end{proof}

\section{Proof of Theorem \ref{t:tslab}} \label{s:between}

In this section, we generalize Proposition \ref{p:pslabn1} to when the
minimal surface is not in a
slab; i.e., we show Theorem \ref{t:tslab}.  $\Sigma^2
\subset B_{c_1 \, r_0}\subset \RR^3$ will be an embedded
minimal disk, $\partial \Sigma \subset
\partial B_{c_1 \, r_0}$, $c_1 \geq 4$, and $y \in \partial B_{2\,
r_0}$.
$\Sigma_1 , \Sigma_2 , \Sigma_3$ will be distinct
components of $B_{r_0}(y) \cap \Sigma$.

\begin{Lem} \label{l:tslab}
Given $\bbeta > 0$,  there exist $2\, c_2 < c_4 < c_3 \leq 1$ so:
Let $\Sigma_3'$ be a component of $B_{ c_3 \, r_0 }(y) \cap
\Sigma_3$ and $y_i \in B_{c_2 \, r_0}(y) \cap \Sigma_i$ for
$i=1,2$.  If $y_1 , y_2$ are in distinct components of $B_{c_4 \,
r_0}(y) \setminus \Sigma_3'$, then there are disjoint stable
embedded minimal disks $\Gamma_{1},\Gamma_2\subset
B_{r_0}(y)\setminus \Sigma$ with $\partial \Gamma_{i} =
\partial  \Sigma_i$,  and (after a rotation)
 graphs $\Gamma_i' \subset \Gamma_i$
over $D_{3 \, c_3 \, r_0}(y)$ so that $y_1 , y_2 , \Sigma_3'$ are
each in their own component of $\Pi^{-1}(D_{3 \, c_3 \, r_0}(y))
\setminus (\Gamma_1' \cup \Gamma_2')$ and $\Gamma_1' , \Gamma_2'
\subset \{ |x_3 - x_3(y)| \leq \bbeta \, c_3 \, r_0 \}$.
\end{Lem}

\begin{proof}
This follows exactly as in the second step of the proof
of Proposition \ref{p:pslabn1}.
\end{proof}

\begin{figure}[htbp]
    \setlength{\captionindent}{20pt}
    \begin{minipage}[t]{0.5\textwidth}
\centering\input{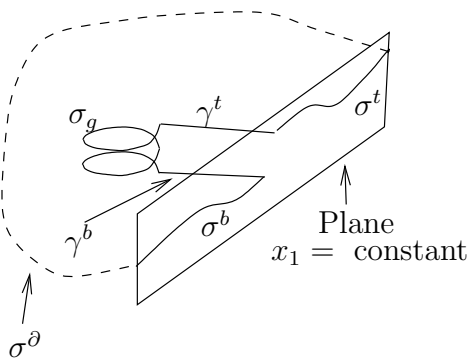}
    \caption{The curve $\gamma_3$
in the proof of Theorem \ref{t:tslab}.
($\gamma_3 = \sigma^b \cup \gamma^b \cup
\sigma_g \cup \gamma^t \cup \sigma^t \cup \sigma^{\partial}$.)}
\label{f:12}    \end{minipage}
\end{figure}

\begin{proof}
(of Theorem \ref{t:tslab}).  Let
$N_1 , \Omega_1, \epsilon > 0$ be given by Theorem
\ref{t:spin4ever} (with $\tau = 1$). Assume that $N_1$ is even.
Let $\beta > 0$ be from
 Proposition \ref{p:cab}.  Set
\begin{equation}    \label{e:bbeta}
    \bbeta = \min \, \{ \beta_s , \, \epsilon , \,
\epsilon /  C_g , \,
    \beta / ( 6 \, [ 60 + 20 \, (N_1 + 3) ] ) \} \, / (5 \, \Omega_1)  \, ,
\end{equation}
where $\beta_s , C_g$ are from Lemma \ref{l:lone}.  Let
 $c_2 , c_3 , c_4$ and
 $\Gamma_i' \subset \Gamma_i$ be given by Lemma \ref{l:tslab}.
Set $c_5  = (60 + 20 \, (N_1+3)) \bbeta \, c_3 / \beta$, so that
$c_5 \leq c_3 / (30 \, \Omega_1)$.  Finally, set $c_1 = 16 \,
\Omega_1$.

We will suppose that $\Sigma_3'$ is not a graph at $z' \in
\Sigma_3'$ and deduce a contradiction.    Set $z = \Pi (z')$.
Since $\Sigma_3'$ separates $y_1,y_2$, it is in the slab between
$\Gamma_1' , \Gamma_2'$.  Using Proposition \ref{p:cab} (with $h=
\bbeta \, c_3 \, r_0 / \beta$) and  \eqr{e:bbeta}, $\Sigma$
contains an $(N_1 + 3)$-valued graph $\Sigma_g$ over $D_{c_3
\,r_0}(z) \setminus D_{c_5 \, r_0}(z)$ and $\Sigma_g$ is also in
the slab. Let $\sigma_g \subset \Sigma_g$ be the $(N_1+2)$-valued
graph over $\partial D_{c_5 \, r_0}(z)$ (see fig. \ref{f:12}).
Let $E$ be the region in
$\Pi^{-1}( D_{c_3 \, r_0/2}(z) \setminus D_{c_3 \, r_0 / (2 \,
\Omega_1) }(z))$ between the sheets of the (concentric)
 $(N_1+1)$-valued subgraph of $\Sigma_g$.

The first step is to find a curve $\gamma_3\subset \Sigma$
containing $\sigma_{g}$ so any stable disk with
boundary $\gamma_3$ is forced to spiral. $\gamma_3$ will have six
pieces: $\sigma_g$, two segments, $\gamma^t , \gamma^b$, in
$\Sigma_g$ which are graphs over a portion of the $\{ x_1 > x_1(z)
\}$ part of the $x_1$-axis, two nodal curves, $\sigma^t ,
\sigma^b$, in $\{ x_1 = {\text{constant}} \}$, and  a
segment $\sigma^{\partial}$ in $\partial \Sigma$. Since $\Sigma_g$
is a graph, there are graphs $\gamma^t , \gamma^b \subset
\Sigma_g$ over a portion of the $\{ x_1 > x_1(z) \}$ part of the
$x_1$-axis from $\partial \sigma_g$ to $y^t , y^b \in \{ x_1 = x_1
(z) + 3\, c_5 \, r_0 \} \cap \Sigma$. By the maximum principle (as
in the proof of Corollary \ref{c:cep}), there are nodal curves
$\sigma^t , \sigma^b \subset \{ x_1 = x_1 (z) + 3\, c_5 \, r_0 \}
\cap \Sigma$ from $y^t , y^b$, respectively, to $y^t_0 , y^b_0 \in
\partial \Sigma$. Finally, connect $y^t_0 , y^b_0$ by a curve
$\sigma^{\partial} \subset
\partial \Sigma$ and set $\gamma_3 = \sigma^b \cup \gamma^b \cup
\sigma_g \cup \gamma^t \cup \sigma^t \cup \sigma^{\partial}$. By
\cite{MeYa}, there is a stable embedded disk $\Gamma
\subset B_{c_1 \, r_0} \setminus \Sigma$ with $\partial \Gamma =
\gamma_3$. Note that $\partial \Gamma \setminus \partial B_{r_0}$
is connected.

We claim that $\sigma^t , \sigma^b$ do not intersect between any
two of the  components  $\{ \sigma_i \}$ of $B_{(c_3 - 2 c_5) \,
r_0 }(z) \cap \{ x_1 = x_1 (z) + 3\, c_5 \, r_0 \} \cap \Sigma_g$.
If not, we can assume that a curve $\sigma \subset \sigma^t$
connects $y^t$ to a point $y_0$ between $\sigma_{i} , \sigma_{i
+1}$. By (a slight variation of) Proposition \ref{p:pslabn1}, the
portion $\Sigma_{y_0}$ of $\Sigma$ between the $i$-th and
$(i+1)$-st sheets of $B_{(c_3 -  c_5) \, r_0}(z) \cap \Sigma_g
\setminus \Pi^{-1} (D_{2 \, c_5 \, r_0}(z))$ is a graph (in fact,
``all the way around''). Note that $B_{3 \, c_5 \, r_0}(z) \cap
\Sigma_{y_0}$ and $B_{3 \, c_5 \, r_0}(z) \cap \Sigma_g$ are in
the same component of $B_{3\, c_5 \, r_0}(z) \cap \Sigma$, since
else the stable disk between them given by \cite{MeYa} would,
using Lemma \ref{l:lone}, intersect $\Sigma_g$. We can therefore
apply the maximum principle as in the proof of Corollary
\ref{c:cep} (i.e., the case $y_0 \in \sigma_j$ for some $j$) to
get the desired contradiction.

We will show next that $\Gamma$ contains an $N_1$-valued graph
$\Gamma_g$ over $D_{c_3  \, r_0/2}(z) \setminus D_{c_3 \, r_0 / (2
\, \Omega_1)}(z)$ with gradient $\leq \epsilon$, $\Pi^{-1}( D_{c_3
\, r_0 / (2 \, \Omega_1)}(z)) \cap (\Gamma_g)^M \subset \{ |x_3 -
x_3(z)| \leq \epsilon \, c_3 \, r_0 / (2 \, \Omega_1) \}$, and a
curve $\eta \subset \Pi^{-1} (D_{c_3 \, r_0 / (2 \, \Omega_1)}(z))
\cap \Gamma \setminus \partial B_{r_0}$ connects $\Gamma_g$ to
$\partial \Gamma \setminus \partial B_{r_0}$. By the previous
paragraph,
\begin{equation}   \label{e:bdfarag}
\dist_{\Gamma} (E \cap \Gamma , \partial \Gamma ) > c_3 \, r_0/(5
\, \Omega_1)  \, .
\end{equation}
 By (the proof of)
Lemma \ref{l:lone} (with $h=c_3 \, r_0/(5 \, \Omega_1)$ and
$\beta= 5 \, \Omega_1 \, \bbeta$),  \eqr{e:bbeta}, and
\eqr{e:bdfarag}, we have that each component of $E \cap \Gamma$ is
a multi-valued graph with gradient $ \leq  5 \, C_g  \, \Omega_1
\, \bbeta \leq \epsilon$. Let $\sigma_c \subset E$ be a graph over
$\partial D_{c_3 \, r_0 / (2 \, \Omega_1)}(z)$.  Using that
$\sigma_c$ separates $\Pi^{-1}(\partial D_{c_3 \, r_0 / (2 \,
\Omega_1)}(z))
  \cap \gamma^t$ and
 $\Pi^{-1}(\partial D_{c_3 \, r_0 / (2 \, \Omega_1)}(z))
  \cap \gamma^b$ in the cylinder $\Pi^{-1}(\partial D_{c_3 \, r_0 / (2 \,
  \Omega_1)}(z))$ (and the description of $\partial \Gamma$),
  there is a   curve $\eta \subset \Pi^{-1} (D_{c_3
\, r_0 / (2 \, \Omega_1)}(z)) \cap \Gamma \setminus \partial
B_{r_0}$ from $\Gamma \cap \sigma_c$ to $\partial \Gamma \setminus
\partial B_{r_0}$.
     Hence, since  $E$ is between the sheets of an
$(N_1+1)$-valued graph, we get the desired $\Gamma_g$.

Combining all of this, Theorem \ref{t:spin4ever} gives
a $2$-valued graph $\Gamma_d \subset \Gamma$  over $D_{c_1 \,
r_0/(2 \, \Omega_1)}(z) \setminus D_{c_3 \, r_0 / (2 \,
\Omega_1)}(z)$ with gradient $\leq 1 $.
Let  $\hat{\gamma}$ be the  component of $B_{(2-2\, c_3) \, r_0}
\cap \gamma$  intersecting $B_{r_0}$.  Note that since
$\partial \gamma = \{ y_1 , y_2\}$ is separated by the slab
between $\Gamma_1' , \Gamma_2'$ and $\gamma \setminus B_{r_0}$ is
$c_2 \, r_0 $-almost-monotone, $\Gamma_d$ separates the endpoints
of $\partial \hat{\gamma}$. Finally, as in the proof of
Proposition \ref{p:pslabn1}, we must have $\Gamma_d \cap
\hat{\gamma} \ne \emptyset$.  This contradiction completes the
proof.
\end{proof}

Many variations of Theorem \ref{t:tslab} hold with
almost the same proof.  One such is:

\begin{Thm} \label{t:tslab2}
There exist $d_1\geq 8$ and $d_2 \leq 1$ so:
 Let $\Sigma^2 \subset B_{d_1 \,r_0}\subset \RR^3$ be an
 embedded minimal disk with $\partial \Sigma
\subset\partial B_{d_1 \, r_0}$ and let $y \in \partial D_{5\,
r_0}$. Suppose that $\Sigma_1 , \Sigma_2 \subset \Sigma$ are disjoint graphs
over $D_{3r_0}(y)$ with gradient $\leq d_2$ and which intersect
$B_{d_2r_0}(y)$.  If they can be connected in $B_{3r_0}\cap
\Sigma$, then any component of $B_{r_0}(y) \cap \Sigma$
which lies between them is a graph.
\end{Thm}

\part{Extending multi-valued graphs off the axis}  \label{p:p4}

In this section $\Sigma \subset
B_{R_0} \subset \RR^3$ will be an embedded
minimal disk with $\partial \Sigma
\subset \partial B_{R_0}$.  In contrast to the results
 of Part \ref{p:p2},
$\Sigma$ is no longer assumed to be stable.

Note that, by \cite{Sc}, we can choose $d_3 >4$ so that: If
$\Gamma_0 \subset B_{d_3 \, s}$ with $\partial \Gamma_0 \subset
\partial B_{d_3 \, s}$ is stable, then each component of $B_{4 \,
s}\cap \Gamma_0$ is a graph (over some plane)
with gradient $\leq 1/2$.

\begin{proof}
(of Theorem \ref{t:spin4ever2}). The proof has two steps.
First, the proof of Theorem \ref{t:tslab} and Lemma \ref{l:harnie}
give a stable disk $\Gamma \subset B_{R_0} \setminus \Sigma$
and a $4$-valued graph
 $\Gamma_4 \subset \Gamma$ so $\Sigma^M$ ``passes between''
$\Gamma_4$. Second, (a slight variation of)
Theorem \ref{t:tslab2}  gives
the  $2$-valued graph $\Sigma_d \subset \Sigma$.

Before proceeding, we choose the constants. Let $C_3 , \delta_2$
be given by Lemma \ref{l:harnie} (with $N=4 $), $d_1 , d_2$ be
from Theorem  \ref{t:tslab2}, and $C_g , \beta_s$ be from Lemma
\ref{l:lone}. Set $\tau_1 = \min \{ \tau / (5 \, C_g) , \beta_s /5
, d_2 / 10 \}$ and $\tau_2 = \min \{ \delta_2 / 3 , \tau_1 / (1 +
3 \, C_3) \}$. Let $N_1 , \Omega_1 , \epsilon$ be given by Theorem
\ref{t:spin4ever} (with $\tau$ there equal to $\tau_2$). For
convenience, assume that $N_1 \geq 16$ is even, $\Omega_1 > 4$,
and rename this $\epsilon$ as $\epsilon_1$.  Set $N= N_1 +3$,
$\Omega = \max \{ d_1 , 8 \, d_3 \, \Omega_1\}$, and
\begin{equation}        \label{e:seteps}
    \epsilon = \min \, \{ \epsilon_1 , \epsilon_1 / (5 \, C_g),
     \beta_s / 5 , 1/4  , d_2 / 10  \}  \, .
\end{equation}
  For $N_2 \leq N$ and $r_0 \leq r_2 < r_3 \leq 1$, 
let $E^{N_2}_{r_2 , r_3}$ be
the region in $\Pi^{-1}(D_{r_3} \setminus D_{r_2
})$ between the sheets of the (concentric)
 $N_2$-valued subgraph of $\Sigma_g$.
Note that $E^{N_2}_{r_2 , r_3}
 \subset \{ x_3^2 \leq \epsilon^2 \, (x_1^2  + x_2^2) \}$.

As in the proof of Theorem \ref{t:tslab}, let $\sigma_g \subset
\Sigma_g$ be an $(N_1+2)$-valued graph over $\partial D_{r_0}$ and
let $\gamma_3\subset \Sigma$  be a curve with six pieces:
$\sigma_g$, two segments, $\gamma^t,\gamma^b$, in $\Sigma_g$ which
are graphs over a portion of the positive part of the $x_1$-axis,
two nodal curves, $\sigma^t , \sigma^b$, in $\{ x_1 = 2\, d_3 \,
r_0 \}$, and $\sigma^{\partial} \subset \partial
\Sigma$. By \cite{MeYa}, there is a stable embedded disk
$\Gamma \subset B_{R_0} \setminus \Sigma$ with $\partial \Gamma =
\gamma_3$.

 Let $\{ \sigma_i \}$ be the  components of $B_{5/8} \cap
\{ x_1 =  2\, d_3 \, r_0 \} \cap \Sigma_g$ and suppose that a
curve $\sigma \subset \sigma^t$ connects $\gamma^t$ to a point
$y_0$ between $\sigma_{i} , \sigma_{i +1}$. By Theorem
\ref{t:tslab2}, the portion $\Sigma_{y_0}$ with $y_0 \in \Sigma_{y_0}$
of $E^{N_1 + 5/2}_{3\,r_0 , 5/8} \cap \Sigma$ is a graph. Note that $B_{d_3 \,
r_0}\cap \Sigma_{y_0}$ and $B_{3 \, r_0}\cap \Sigma_g$ are in the
same component of $B_{d_3 \, r_0}\cap \Sigma$, since else the
stable disk between them given by \cite{MeYa} would intersect
$\Sigma_g$ (using \cite{Sc}). Applying the maximum principle as
before gives the desired contradiction. Hence, $\sigma^t ,
\sigma^b$ do not intersect between any of the $\sigma_i$'s.
Therefore, if $z \in E^{N_1 +1}_{4d_3 \, r_0,1/2} \cap \Gamma$, then
\begin{equation}   \label{e:bdfarag2}
\dist_{\Gamma} (z , \partial \Gamma ) \geq |\Pi(z)| / 4 \, .
\end{equation}
   By the same linking argument as before, $E^{N_1 +1}_{4d_3 \, r_0,1/2}  \cap
\Gamma$ contains an $N_1$-valued graph $\Gamma_g$ over $D_{1/2}
\setminus D_{4 \, d_3 \, r_0 }$ with gradient $\leq 5 \, C_g \,
\epsilon$, $\Pi^{-1}( \partial D_{ 4 \, d_3 \, r_0 }) \cap
\Gamma_g \subset \{ |x_3| \leq 4\, \epsilon \, d_3 \, r_0  \}$,
and a curve $\eta \subset \Pi^{-1} (D_{4 \, d_3 \, r_0 }) \cap
\Gamma \setminus \partial B_{R_0}$ connects $\Gamma_g$ to
$\partial \Gamma \setminus \partial B_{R_0}$.
 Since $\Omega_1 < 1/ (8 \, d_3 \, r_0)$,
 Theorem \ref{t:spin4ever} implies that
$\Gamma$ contains a $2$-valued graph $\Gamma_d$ over $D_{
R_0/\Omega_1} \setminus D_{4 \, d_3 \, r_0}$ with gradient $\leq
\tau_2 < 1$. In particular, $\Gamma_d \subset \{ x_3^2 \leq
\tau_2^2 ( x_1^2 + x_2^2) \}$.
  Next, we apply  Lemma \ref{l:harnie}
 to extend $\Gamma_d$ to a
$4$-valued graph $\Gamma_4$ over $D_{ 5\, R_0/(6 \, \Omega_1)}
\setminus D_{5 \, d_3 \, r_0}$ with gradient
$\leq \tau_2 + 3 \, C_3 \,  \tau_2 \leq \tau_1$.
 Let
$E_{\Gamma}$ be the region in $\Pi^{-1}(D_{ R_0/(2 \, \Omega_1)}
\setminus D_{15 \, d_3 \, r_0})$ between the sheets of the
(concentric) $3$-valued subgraph
 of $\Gamma_4$, so that
$E_{\Gamma} \subset \{ x_3^2 \leq \tau_1^2 \, (x_1^2  + x_2^2) \}$.

  If $z \in E_{\Gamma} \cap \Sigma$, then
there is a curve $\gamma_z \subset \Gamma_4$ with each component
of $\gamma_z \setminus \Pi^{-1}(D_{5 \, d_3 \, r_0})$ a graph over the segment
$\gamma_{0,z}$, $\partial \gamma_z = \{ y_z^2 , y_z^4 \}$, and
$y_z^2 , y_z^4$ are in distinct components of
$B_{3 \, |\Pi(z)|/ 5}(\Pi(z)) \cap \Gamma$ with $z$
between these components.
 By (a slight variation of) Theorem \ref{t:tslab2}
(using $\Sigma \cup \Gamma$ as a barrier rather than just
$\Sigma$), the portion of $\Sigma$ inside $B_{R_0/d_1} \cap
E_{\Gamma}$ is a graph over $\Gamma_4$.  This is nonempty since
$(\Sigma_g)^M$ begins in $E_{\Gamma}$, so we get the desired
$2$-valued graph $\Sigma_d$ with gradient $\leq 5 \, C_g \, \tau_1
\leq \tau$ (by
 Lemma \ref{l:lone}).
\end{proof}

\appendix

\section{Catenoid foliations} \label{s:cat}

We recall here some consequences of the maximum principle for an
embedded minimal surface $\Sigma$ in a slab. Let $\cat (y)$ be the
vertical catenoid centered at $y=(y_1 , y_2 , y_3)$ given by
\begin{equation}
\cat (y) = \{ x \in \RR^3 \, | \, \cosh^2 (x_3 - y_3) = (x_1-
y_1)^2 + (x_2-y_2)^2 \} \, .
\end{equation}
Given $0 < \theta < \pi/2$, let $\partial N_{\theta} (y)$ be the
cone
\begin{equation}
\{ x \, | \, (x_3 - y_3)^2 = |x-y|^2 \, \sin^2 \theta \} \, .
\end{equation}
Since $\cosh t > t$ for $t\geq 0$, it follows that $\partial
N_{\pi / 4}(y) \cap \cat(y) = \emptyset$. Set
\begin{equation} \label{e:deft0}
\theta_0 = \inf \, \{ \theta \, | \, \partial N_{\theta}(y) \cap
\cat(y) = \emptyset \} \, ,
\end{equation}
so that $ \partial N_{\theta_0}(y)$ and $\cat (y)$ intersect
tangentially in a pair of circles. Let $\cat_0(y)$ be the
component of $\cat(y) \setminus \partial N_{\theta_0}(y)$
containing the neck $\{ x \, | \, x_3 = y_3 , (x_1- y_1)^2 +
(x_2-y_2)^2 = 1 \}$. If $x \in \cat_0(y)$, then $\overline{y,x}
\cap \cat_0(y) = \{ x \}$ since $\cosh$ is convex and $\cosh'(0) =
0$; i.e., $\cat_0(y)$ is a radial graph. Dilating $\cat_0 (y)$ (see
fig. \ref{f:shn12})
gives a minimal foliation of the solid (open) cone
\begin{equation}
N_{\theta_0} (y) = \{ x \, | \, (x_3 - y_3)^2 < |x-y|^2 \, \sin^2
\theta_0 \} \, .
\end{equation}
The leaves have boundary in $\partial N_{\theta_0} (y)$ and are
level sets of the function $f_y$ given by
\begin{equation} \label{e:deffy}
y + (x - y) / f_y (x) \in \cat_0 (y) \, .
\end{equation}
Choose $\beta_A > 0$ small so that $\{ x \, | \, |x_3 -y_3| \leq
2\, \beta_A \, h \} \setminus B_{h / 8}(y) \subset N_{\theta_0}
(y)$ and
\begin{equation}
\{ x \, | \, f_y (x) = 3 \, h / 16 \} \cap \{x \, | \, |x_3 - y_3|
\leq 2 \, \beta_A \, h \} \subset B_{ 7 \, h / 32}(y) \, .
\label{e:lowh}
\end{equation}

\begin{figure}[htbp]
    \setlength{\captionindent}{20pt}
    \begin{minipage}[t]{0.5\textwidth}
    \centering\input{shn12.pstex_t}
    \caption{The catenoid foliation.}
    \label{f:shn12}    \end{minipage}\begin{minipage}[t]{0.5\textwidth}
\centering\input{shn13.pstex_t}
    \caption{An $n$-prong singularity.}
\label{f:shn13}    \end{minipage}
\end{figure}

The intersection of two embedded minimal surfaces is locally given
by $2n$ embedded arcs meeting at equal angles as in fig. \ref{f:shn13},
i.e., an
``$n$-prong singularity'' (e.g., the set where $(x+ i y)^n$ is
real); see claim $1$ in lemma $4$ of \cite{HoMe}. This immediately
implies:

\begin{Lem} \label{l:apa}
If $z \in \Sigma \subset N_{\theta_0} (y)$ is a nontrivial
interior critical point of $f_y |_{\Sigma}$, then $\{ x \in \Sigma
\, | \, f_y (x) = f_y (z) \}$ has an $n$-prong singularity at $z$
with $n \geq 2$.
\end{Lem}

As a consequence, we get a version of the usual strong maximum
principle:

\begin{Lem} \label{l:apb}
If $\Sigma \subset N_{\theta_0} (y)$, then $f_y |_{\Sigma}$ has no
nontrivial interior local extrema.
\end{Lem}

\begin{Cor} \label{c:c1}
If $\Sigma \subset B_h(y) \cap \{ x \, | \, |x_3 - y_3| \leq 2\,
\beta_A \, h \}$, $\partial \Sigma \subset \partial
B_h(y)$, and $ B_{3 \, h / 4}(y) \cap \Sigma \ne \emptyset$, then
$B_{ h / 4}(y) \cap \Sigma \ne \emptyset $.
\end{Cor}

\begin{proof}
Scaling \eqr{e:lowh} by $4$, $\{ x \in \Sigma \, | \, f_y (x) = 3
\, h / 4 \} \subset B_{7 \, h / 8}(y) \setminus B_{3 \, h /
4}(y)$. By Lemma \ref{l:apb}, $f_y$ has no interior minima in
$\Sigma$ so the corollary now follows from $f_y(x) \leq |x-y|$.
\end{proof}

Iterating Corollary \ref{c:c1} along a chain of balls gives:

\begin{Cor} \label{c:c3}
If $\Sigma \subset \{ |x_3| \leq 2\, \beta_A \, h \}$,
$p,q \in \{ x_3 = 0\}$, $T_h (\gamma_{p,q}) \cap \partial \Sigma =
\emptyset$, and $y_p \in B_{h/4}(p) \cap \Sigma$, then a curve
$\cag \subset T_{h} (\gamma_{p,q}) \cap \Sigma$ connects $y_p$ to
$B_{h / 4}(q) \cap \Sigma$.
\end{Cor}

\begin{proof}
Choose $y_0 = p , y_1 , y_2 , \dots , y_n = q \in \gamma_{p , q}$
with $|y_{i-1} - y_{i}| = h / 2$ for $ i < n$ and $|y_{n-1} -
y_{n}| \leq h / 2$. Repeatedly applying Corollary \ref{c:c1} for
$1 \leq i \leq n$, gives $ \cag_i : [0,1] \to B_{ h }(y_i) \cap
\Sigma$ with $\cag_1 (0) = y_p$, $\cag_i (1) \in B_{ h / 4}(y_{i})
\cap \Sigma$, and $\cag_{i+1} (0) = \cag_i (1)$. Set $\cag =
\cup_{i=1}^{n} \, \cag_i$.
\end{proof}

This produces curves which are ``$h$-almost monotone'' in the
sense that if $y \in \cag$, then $B_{4 \, h}(y) \cap \cag$ has
only one component which intersects $B_{2\, h }(y)$.

\begin{Cor} \label{c:c2}
If $\Sigma \subset \{ |x_3| \leq 2 \, \beta_A \, h \}$
and $E$ is an unbounded component of $\RR^2 \setminus T_{h / 4}
(\Pi (\partial \Sigma))$, then $ \Pi (\Sigma) \cap E = \emptyset$.
\end{Cor}

\begin{proof}
Given $y \in E$, choose a curve $\gamma : [0, 1] \to \RR^2
\setminus T_{h / 4} (\Pi (\partial \Sigma))$ with $| \gamma (0) |
> \sup_{x \in \Sigma} |x| + h$ and $\gamma (1) = y$. Set $\Sigma_t
= \{x \in \Sigma \, | \, f_{\gamma (t)}(x) = 3 \, h / 16 \}$. By
\eqr{e:lowh}, $\Sigma_t \subset B_{7 \, h / 32}(\gamma(t)) $, so
that $\Sigma_0 = \emptyset$ and $\Sigma_t \cap \partial \Sigma =
\emptyset$. By Lemma \ref{l:apb}, either $\Sigma_t = \emptyset$ or
$\Sigma_t$ contains an arc of transverse intersection. In
particular, there cannot be a first $t > 0$ with $\Sigma_t \ne
\emptyset$, giving the corollary.
\end{proof}

\end{document}